%% file: EMBiventricular_CMAME.tex
\documentclass[preprint, 3p]{elsarticle}

\input{packages.tex}
\input{macros.tex}

\graphicspath{{pictures/}}
\doublefalse

\begin{document}
	
	\begin{frontmatter}
		\title{{\papertitle}}
		
		\input{header_CMAME.tex}
		
		\begin{abstract}
			\input{parts_abstract.tex}
		\end{abstract}
		
		
	\end{frontmatter}
	
	\input{parts_intro.tex}
	\input{parts_models.tex}
	\input{parts_numerics.tex}
	\input{parts_results.tex}
	\input{parts_conclusions.tex}

	\input{acknowledgements.tex}
	
	\begin{appendices}
		\input{parts_app_params.tex}
		\input{parts_app_BC.tex}
		\input{parts_app_schur.tex}
	\end{appendices}
	
	\bibliographystyle{elsarticle-num}
	\bibliography{references}
	
\end{document}

%% file: packages.tex
\usepackage[english]{babel}
\usepackage[utf8]{inputenc}
\usepackage[colorinlistoftodos]{todonotes}
\usepackage{amsmath}
\usepackage{amssymb}
\usepackage{amsthm} 
\usepackage{mathtools}
\usepackage{mathrsfs}
\usepackage{siunitx}
\usepackage{gensymb}
\usepackage{comment}

\usepackage[hidelinks]{hyperref}
\usepackage{algorithm}
\usepackage{algpseudocode}

\usepackage{subfig}

\usepackage{tikz}

\usepackage{tabularx}
\usepackage{booktabs}

\usepackage{empheq}

\usepackage[toc,page]{appendix}

\theoremstyle{plain}
\theoremstyle{plain}
\theoremstyle{plain}
\theoremstyle{plain}
\theoremstyle{plain}
\theoremstyle{definition}
\theoremstyle{remark}

\everymath{\displaystyle}

%% file: macros.tex
\newif\ifdouble

\newcommand{\papertitle}{
3D-0D closed-loop model for the simulation of cardiac biventricular electromechanics
}

\newcommand{\keywordOne}{Cardiac electromechanics}
\newcommand{\keywordTwo}{Cardiac fiber architecture}
\newcommand{\keywordThree}{Multiphysics modeling}
\newcommand{\keywordFour}{Finite Elements}
\newcommand{\keywordFive}{3D-0D coupling}

\newcommand{\hone}{\mathrm{h}_\mathrm{1}}
\newcommand{\htwo}{\mathrm{h}_\mathrm{2}}

\newcommand{\EPnumGating}{n_{\boldsymbol{w}}}

\newcommand{\ActNumVariables}{n_{\mathbf{s}}}
\newcommand{\CircNumVariables}{n_{\mathbf{c}}}

\newcommand{\PhyEP}{\mathscr{E}}

\newcommand{\PhyEPIMEX}{\mathscr{E}_{\mathrm{IMEX}}}

\newcommand{\PhyMec}{\mathscr{M}}
\newcommand{\PhyMecI}{\mathscr{M}_{\mathrm{I}}}
\newcommand{\PhyAct}{\mathscr{A}}
\newcommand{\PhyActE}{\mathscr{A}_{\mathrm{E}}}
\newcommand{\PhyCirc}{\mathscr{C}}
\newcommand{\PhyCircE}{\mathscr{C}_{\mathrm{E}}}
\newcommand{\PhyCoupl}{\mathscr{V}} 
\newcommand{\PhyCouplI}{\mathscr{V}_{\mathrm{I}}} 

\newcommand{\VLVthreedim}{V_{\mathrm{LV}}^{\mathrm{3D}}}

\newcommand{\VRVthreedim}{V_{\mathrm{RV}}^{\mathrm{3D}}}
\newcommand{\Vthreedimi}{V_{\mathrm{i}}^{\mathrm{3D}}}

\newcommand{\VLA}{V_{\mathrm{LA}}}
\newcommand{\VLV}{V_{\mathrm{LV}}}
\newcommand{\VRA}{V_{\mathrm{RA}}}
\newcommand{\VRV}{V_{\mathrm{RV}}}

\newcommand{\VnLA}{V_{\mathrm{0,LA}}}

\newcommand{\VnRA}{V_{\mathrm{0,RA}}}

\newcommand{\PLV}{p_{\mathrm{LV}}}

\newcommand{\PRV}{p_{\mathrm{RV}}}

\newcommand{\EpLA}{E_{\mathrm{LA}}^{\mathrm{A}}}

\newcommand{\EpRA}{E_{\mathrm{RA}}^{\mathrm{A}}}

\newcommand{\EaLA}{E_{\mathrm{LA}}^{\mathrm{B}}}

\newcommand{\EaRA}{E_{\mathrm{RA}}^{\mathrm{B}}}

\newcommand{\TacLA}{T_{\mathrm{LA}}^{\mathrm{ac}}}
\newcommand{\TacRA}{T_{\mathrm{RA}}^{\mathrm{ac}}}
\newcommand{\tacLA}{t_{\mathrm{LA}}^{\mathrm{ac}}}
\newcommand{\tacRA}{t_{\mathrm{RA}}^{\mathrm{ac}}}
\newcommand{\TarLA}{T_{\mathrm{LA}}^{\mathrm{ar}}}
\newcommand{\TarRA}{T_{\mathrm{RA}}^{\mathrm{ar}}}

\newcommand{\QarSYS}{Q_{\mathrm{AR}}^{\mathrm{SYS}}}
\newcommand{\QarPUL}{Q_{\mathrm{AR}}^{\mathrm{PUL}}}
\newcommand{\QvnSYS}{Q_{\mathrm{VEN}}^{\mathrm{SYS}}}
\newcommand{\QvnPUL}{Q_{\mathrm{VEN}}^{\mathrm{PUL}}}

\newcommand{\CarSYS}{C_{\mathrm{AR}}^{\mathrm{SYS}}}
\newcommand{\CarPUL}{C_{\mathrm{AR}}^{\mathrm{PUL}}}
\newcommand{\CvnSYS}{C_{\mathrm{VEN}}^{\mathrm{SYS}}}
\newcommand{\CvnPUL}{C_{\mathrm{VEN}}^{\mathrm{PUL}}}

\newcommand{\ParSYS}{p_{\mathrm{AR}}^{\mathrm{SYS}}}
\newcommand{\ParPUL}{p_{\mathrm{AR}}^{\mathrm{PUL}}}
\newcommand{\PvnSYS}{p_{\mathrm{VEN}}^{\mathrm{SYS}}}
\newcommand{\PvnPUL}{p_{\mathrm{VEN}}^{\mathrm{PUL}}}

\newcommand{\RarSYS}{R_{\mathrm{AR}}^{\mathrm{SYS}}}
\newcommand{\RarPUL}{R_{\mathrm{AR}}^{\mathrm{PUL}}}
\newcommand{\RvnSYS}{R_{\mathrm{VEN}}^{\mathrm{SYS}}}
\newcommand{\RvnPUL}{R_{\mathrm{VEN}}^{\mathrm{PUL}}}

\newcommand{\LarSYS}{L_{\mathrm{AR}}^{\mathrm{SYS}}}
\newcommand{\LarPUL}{L_{\mathrm{AR}}^{\mathrm{PUL}}}
\newcommand{\LvnSYS}{L_{\mathrm{VEN}}^{\mathrm{SYS}}}
\newcommand{\LvnPUL}{L_{\mathrm{VEN}}^{\mathrm{PUL}}}

\newcommand{\Rmin}{R_{\mathrm{min}}}
\newcommand{\Rmax}{R_{\mathrm{max}}}

\newcommand{\Circ}{\boldsymbol{c}}

\newcommand{\CircInit}{\boldsymbol{c}_{0}}

\newcommand{\CircRhs}{\boldsymbol{D}}

\newcommand{\fZero}{{\mathbf{f}_0}}
\newcommand{\sZero}{{\mathbf{s}_0}}
\newcommand{\nZero}{{\mathbf{n}_0}}

\newcommand{\GammaA}{\Gamma_0^{\mathrm{a}}}
\newcommand{\GammaB}{\Gamma_0^{\mathrm{b}}}
\newcommand{\GammaN}{\Gamma_0^{\mathrm{n}}}

\newcommand{\GammaClosedA}{\overline{\Gamma}_0^{\mathrm{a}}}
\newcommand{\GammaClosedB}{\overline{\Gamma}_0^{\mathrm{b}}}
\newcommand{\GammaClosedN}{\overline{\Gamma}_0^{\mathrm{n}}}

\newcommand{\chiA}{\chi_{\mathrm{a}}}

\newcommand{\chiB}{\chi_{\mathrm{b}}}

\newcommand{\XiHat}{\hat{\xi}}

\newcommand{\eL}{\widehat{\boldsymbol{e}}_{\mathrm{\ell}}}
\newcommand{\eT}{\widehat{\boldsymbol{e}}_{\mathrm{t}}}
\newcommand{\eN}{\widehat{\boldsymbol{e}}_{\mathrm{n}}}

\newcommand{\ii}{\mathrm{i}}

\newcommand{\alphai}{\alpha_{\mathrm{i}}}
\newcommand{\betai}{\beta_{\mathrm{i}}}
\newcommand{\thetaEndoi}{\theta_{\mathrm{endo,i}}}
\newcommand{\thetaEpii}{\theta_{\mathrm{epi,i}}}

\newcommand{\EPchim}{\chi_\mathrm{m}}
\newcommand{\EPCm}{C_\mathrm{m}}
\newcommand{\EPIion}{{\mathcal{I}_{\mathrm{ion}}}}
\newcommand{\EPIapp}{{\mathcal{I}_{\mathrm{app}}}}
\newcommand{\tIapp}{{t_{\mathrm{app}}}}
\newcommand{\tIappzLV}{{t^0_{\mathrm{LV,app}}}}
\newcommand{\tIappzRV}{{t^0_{\mathrm{RV,app}}}}
\newcommand{\DIapp}{{\delta_{\mathrm{app}}}}

\newcommand{\EPdiffTens}{\boldsymbol{D}}

\newcommand{\EPrhsGating}{\boldsymbol{H}}

\newcommand{\LinPlus}{{\mathrm{Lin}^+}}
\newcommand{\Inv}[1]{{\mathcal{I}_{#1}}}
\newcommand{\IIVf}{\Inv{4f}}

\newcommand{\IIVs}{\Inv{4s}}
\newcommand{\IIVn}{\Inv{4n}}

\newcommand{\displ}{\mathbf{d}} 

\newcommand{\BCmecCepiN}{{C_\bot^{\mathrm{epi}}}}
\newcommand{\BCmecKepiN}{{K_\bot^{\mathrm{epi}}}}
\newcommand{\BCmecCepiT}{{C_\parallel^{\mathrm{epi}}}}
\newcommand{\BCmecKepiT}{{K_\parallel^{\mathrm{epi}}}}

\newcommand{\BCmecKepiTens}{\mathbf{K}^{\mathrm{epi}}}
\newcommand{\BCmecCepiTens}{\mathbf{C}^{\mathrm{epi}}}

\newcommand{\GammaBase}{\Gamma_0^{\mathrm{base}}}
\newcommand{\GammaEpi}{\Gamma_0^{\mathrm{epi}}}
\newcommand{\GammaEndo}{\Gamma_0^{\mathrm{endo}}}
\newcommand{\GammaEndoi}{\Gamma_0^{\mathrm{endo,i}}}
\newcommand{\GammaEndoRV}{\Gamma_0^{\mathrm{endo,RV}}}
\newcommand{\GammaEndoLV}{\Gamma_0^{\mathrm{endo,LV}}}
\newcommand{\GammaBaset}{\Gamma_t^{\mathrm{base}}}

\newcommand{\GammaClosedBase}{\overline{\Gamma}_0^{\mathrm{base}}}
\newcommand{\GammaClosedEpi}{\overline{\Gamma}_0^{\mathrm{epi}}}
\newcommand{\GammaClosedEndoRV}{\overline{\Gamma}_0^{\mathrm{endo,rv}}}
\newcommand{\GammaClosedEndoLV}{\overline{\Gamma}_0^{\mathrm{endo,lv}}}

\newcommand{\mecF}{\mathbf{F}}
\newcommand{\mecC}{\mathbf{C}}
\newcommand{\tenspiola}{\mathbf{P}}
\newcommand{\tenscaychy}{\mathbf{T}}
\newcommand{\identity}{\mathbf{I}}

\newcommand{\mecNref}{{\mathbf{N}}}
\newcommand{\mecNact}{{\mathbf{n}}}

\newcommand{\stress}{S}

\newcommand{\BCmecVbaseLV}{\mathbf{v}_{\mathrm{LV}}^{\mathrm{base}}}
\newcommand{\BCmecVbaseRV}{\mathbf{v}_{\mathrm{RV}}^{\mathrm{base}}}
\newcommand{\BCmecVbaseDiscr}[1]{\mathbf{v}^{\mathrm{base,n}}_{#1}}

\newcommand{\Cr}{C_{\mathrm{lrv}}}

\newcommand{\BCweight}{\phi}

\newcommand{\GammaBasetTilde}{\widetilde{\Gamma}_t^{\text{base}}}
\newcommand{\GammaEpitTilde}{\widetilde{\Gamma}_t^{\text{epi}}}

\newcommand{\GammaBaseL}{\Gamma_0^{\text{base,L}}}
\newcommand{\GammaEndoL}{\Gamma_0^{\text{endo,L}}}

\newcommand{\GammaEndotL}{\Gamma_t^{\text{endo,L}}}

\newcommand{\OmegaFluidtL}{\Omega_t^{\text{fluid,L}}}
\newcommand{\OmegaFluidtTildeL}{\widetilde{\Omega}_t^{\text{fluid,L}}}

\newcommand{\GammaEndotTildeL}{\widetilde{\Gamma}_t^{\text{endo,L}}}

\newcommand{\GammaBaseR}{\Gamma_0^{\text{base,R}}}
\newcommand{\GammaEndoR}{\Gamma_0^{\text{endo,R}}}

\newcommand{\GammaEndotR}{\Gamma_t^{\text{endo,R}}}

\newcommand{\OmegaFluidtR}{\Omega_t^{\text{fluid,R}}}
\newcommand{\OmegaFluidtTildeR}{\widetilde{\Omega}_t^{\text{fluid,R}}}

\newcommand{\GammaEndotTildeR}{\widetilde{\Gamma}_t^{\text{endo,R}}}

\newcommand{\adofs}{\mathbf{a}}
\newcommand{\displdofs}{\displ}

\newcommand{\Cai}{{[\mathrm{Ca}^{2+}]_{\mathrm{i}}}}

\newcommand{\SL}{{SL}} 

\newcommand{\Tens}{T_\mathrm{a}} 
\newcommand{\TensSubscriptDofs}[1]{\mathbf{T}_{\mathrm{a},{#1}}} 
\newcommand{\TensMax}{T_\mathrm{a}^{\mathrm{max}}} 
\newcommand{\ActStateHF}{\mathbf{s}}

\newcommand{\ActRhs}{\boldsymbol{K}}

\newcommand{\Nsub}{N_{\mathrm{sub}}}

\newcommand{\OmegaRef}{\widetilde{\Omega}}

\newcommand{\pressLVRef}{\widetilde{p}_{\mathrm{LV}}}
\newcommand{\pressRVRef}{\widetilde{p}_{\mathrm{RV}}}
\newcommand{\tensRef}{\widetilde{T}_{\mathrm{a}}}

%% file: header_CMAME.tex
\journal{Computer Methods in Applied Mechanics and Engineering }

\author[mox]{Roberto~Piersanti}
\author[mox]{Francesco~Regazzoni\corref{cor1}}\ead{francesco.regazzoni@polimi.it}
\author[mox]{Matteo~Salvador}
\author[usa]{Antonio~F.~Corno}
\author[mox]{Luca~Dede'} 
\author[labs]{Christian~Vergara}
\author[mox,epfl]{Alfio~Quarteroni}

\address[mox]{MOX - Dipartimento di Matematica, Politecnico di Milano,\par P.zza Leonardo da Vinci 32, 20133 Milano, Italy}
\address[usa]{Houston Children's Heart Institute, Hermann Children's Hospital \par University of Texas Health, McGovern Medical School, Houston, Texas, USA}
\address[labs]{LaBS, Dipartimento di Chimica, Materiali e Ingegneria Chimica "Giulio Natta", Politecnico di Milano,\par P.zza Leonardo da Vinci 32, 20133 Milano, Italy}
\address[epfl]{Mathematics Institute, \'{E}cole Polytechnique F\'{e}d\'{e}rale de Lausanne,\par Av. Piccard, CH-1015 Lausanne, Switzerland (\textit{Professor Emeritus})}

\cortext[cor1]{Corresponding author.}

%% file: parts_abstract.tex
Two crucial factors for accurate numerical simulations of cardiac electromechanics,
which are also essential to reproduce the synchronous activity of the heart, are: i) accounting for the interaction between the heart and the circulatory system that determines pressures and volumes loads in the heart chambers; ii)~reconstructing the muscular fiber architecture that drives the electrophysiology signal and the myocardium contraction. 
In this work, we present a 3D biventricular electromechanical model coupled with a 0D closed-loop model of the whole cardiovascular system that addresses the two former crucial factors. With this aim, we introduce a boundary condition for the mechanical problem that accounts for the neglected part of the domain located on top of the biventricular basal plane and that is consistent with the principles of momentum and energy conservation.
We also discuss in detail the coupling conditions that stand behind the 3D and the 0D models.
We perform electromechanical simulations in physiological conditions using the 3D-0D model and we show that our results match the experimental data of relevant mechanical biomarkers available in literature. Furthermore, we investigate different arrangements in cross-fibers active contraction.
We prove that an active tension along the sheet direction counteracts the myofiber contraction, while the one along the sheet-normal direction enhances the cardiac work. 
Finally, several myofiber architectures are analysed. We show that a different fiber field in the septal area and in the transmural wall effect the pumping functionality of the left ventricle. 

%% file: parts_intro.tex
\section{Introduction}
\label{sec: introduction}

Over the years, computational models of cardiac electromechanics (EM) ~\cite{gurev2011models,augustin2016anatomically,augustin2020impact,land2018influence,nordsletten2011coupling,strocchi2020publicly,math9111247} have been developed with increasingly biophysical detail, by taking into account the interacting physical phenomena 
characteristic of the heart EM - electrophysiology, active contraction, mechanics~\cite{sermesant2012patient,peirlinck2021precision,smith2004multiscale,chabiniok2016multiphysics,crampin2004computational}.
However, most of the existing EM models refer to the left ventricle (LV) only~\cite{marx2020personalization,Gerbi,regazzoni2020model,salvador2020intergrid,levrero2020sensitivity,propp2020orthotropic,rossi2014thermodynamically} and neglect the important effects of the right ventricular deformation on the heart pumping function \cite{palit2015computational}. Only recently, biventricular EM models~\cite{sermesant2005simulation,chapelle2009numerical,goktepe2010electromechanics,crozier2016image,hirschvogel2017monolithic,ahmad2018multiphysics,garcia2019towards,augustin2020physiologically} have been purposely developed.
Two crucial factors for an accurate numerical simulation of the cardiac EM, which are also essential to reproduce the synchronous activity of the heart, are: i) accounting for the interaction between the heart and the circulatory system and ii) reconstructing the muscular fiber architecture.

The coupling between the circulatory system haemodynamics and the cardiac mechanics determines pressures and
volumes in the heart chambers~\cite{Gerbi,augustin2020physiologically,kerckhoffs2007coupling,dede2021modeling,guan2020effect,wang2021human}. Typically, 3D EM models are coupled with Windkessel-type preload/afterload models for the circulatory system~\cite{elzinga1973pressure,liu2015multi,segers2008three,stergiopulos1999total,wang2003time,westerhof1991normalized}.
In these models, the different phases of the pressure-volume loop (PV-loop) are managed by solving different sets of differential  equations, one for each phase~\cite{Gerbi,dede2020segregated,eriksson2013influence,usyk2002computational}.
Still, more meaningful and physiologically sound interface conditions can be obtained by coupling the 3D EM model with a 0D closed-loop fluid dynamics model of the complete circulatory system for the whole cardiac cycle.~\cite{regazzoni2020model,blanco20103d,arts2005adaptation,neal2007subject,paeme2011mathematical}. A further advantage of the latter approach is that closed-loop circulation models do not require to be adapted through the different phases of the cardiac cycle~\cite{nordsletten2011coupling,augustin2020physiologically,dede2021modeling,regazzoni2020numerical}. 
However, solving efficiently the coupling between the EM model and the closed-loop model for the whole cardiovascular system is a challenging task~\cite{augustin2020physiologically}.
To the best of our knowledge, this coupled problem has been so far addressed only in a few works, namely~\cite{regazzoni2020model,hirschvogel2017monolithic,augustin2020physiologically,kerckhoffs2007coupling}.

The myocardial fibers plays a key role in the electric signal propagation and in the myocardial contraction~\cite{piersanti2021modeling,streeter1969fiber,roberts1979influence,gil2019influence,carreras2012left}. Due to the difficulty of reconstructing cardiac fibers from medical imaging, a widely used strategy for generating myofiber orientations in EM models relies on the so called Laplace–Dirichlet-Rule-Based-Methods (LDRBMs)~\cite{bayer2012novel,wong2014generating,doste2019rule,quarteroni2017integrated}, recently analysed under a communal mathematical setting~\cite{piersanti2021modeling}.
While it is well recognized that myofibers orientation is crucial for the construction of a realistic EM model, their architecture has been explored only in a few works and it is not fully understood~\cite{palit2015computational,guan2020effect,gil2019influence,azzolin2020effect,pluijmert2017determinants}. 
Another crucial issue for the reconstruction of a suitable cardiac fiber architecture
consists in considering the myofibers dispersion around a predominant direction~\cite{guan2020effect,guan2021modelling,ahmad2018region,sommer2015biomechanical}. Based on experimental measures \cite{lin1998multiaxial}, cross-fibers active tension has been introduced in~\cite{genet2014distribution,sack2018construction,wenk2012first} to model the contraction caused by dispersed myofibers. However, to the best of our knowledge, this aspect was addressed in EM models only in~\cite{levrero2020sensitivity,eriksson2013modeling}. 

With the aim of facing the computational challenges formerly described, our contributions in this paper move along two strands: 
i) on the one hand, we present a biophysically detailed 3D biventricular EM model coupled with a 0D closed-loop lumped parameters model for the haemodynamics of the whole circulatory system; ii) on the other hand, we investigate the effect of different myofiber architectures, by considering three type of LDRBMs, on the biventricular EM. Specifically, we provide the mathematical formulation and the numerical framework of the coupled 3D-0D model carefully inspecting the coupling conditions of these heterogeneous models. We propose an effective boundary condition for the
mechanical problem that accounts for the neglected part of the domain located above the biventricular basal plane and that fulfils the principles of momentum and energy conservation. We report the results of several electromechanical simulations in physiological conditions using the proposed 3D-0D model. Our results match the experimental data of relevant mechanical biomarkers available in the literature~\cite{maceira2006normalized,tamborini2010reference,maceira2006reference,sugimoto2017echocardiographic,bishop1997clinical,emilsson2006mitral,sechtem1987regional}.  
Furthermore, we study at which extent different configurations in cross-fibers active contraction, that surrogate the myofibers dispersion, affect the electromechanical simulations.

This paper is organized as follows.
In Section~\ref{sec: mathematicalmodeling} we briefly recall the fiber generation methods used to model the cardiac muscle fiber architecture in biventricular geometries.
Moreover, we fully present the mathematical formulation for the closed-loop 3D-0D EM model.
Then, in Section~\ref{sec: numericaldiscretization} we present the numerical approximation of the 3D-0D model along with the coupling strategy.
In Section~\ref{sec: numericalresults}, we show the numerical results obtained with the proposed model.
Finally, in Section~\ref{sec: conclusions} we draw our conclusions.

%% file: parts_models.tex
\section{Mathematical models}
\label{sec: mathematicalmodeling}
In this section we provide a brief overview of the fiber generation methods used to reconstruct the cardiac muscle architecture in biventricular geometries (Section~\ref{subsec: fibers}) and we fully present the 3D cardiac EM model for the human heart function together with a 0D model of the whole cardiovascular system~(Section~\ref{subsec: elettromech_model}). 
Finally, we show the strategy to reconstruct the unloaded (i.e. stress-free) configuration~(Section~\ref{subsec: refconf}).

\begin{figure}[t!]
	\centering
	\includegraphics[keepaspectratio, width=0.85\textwidth]{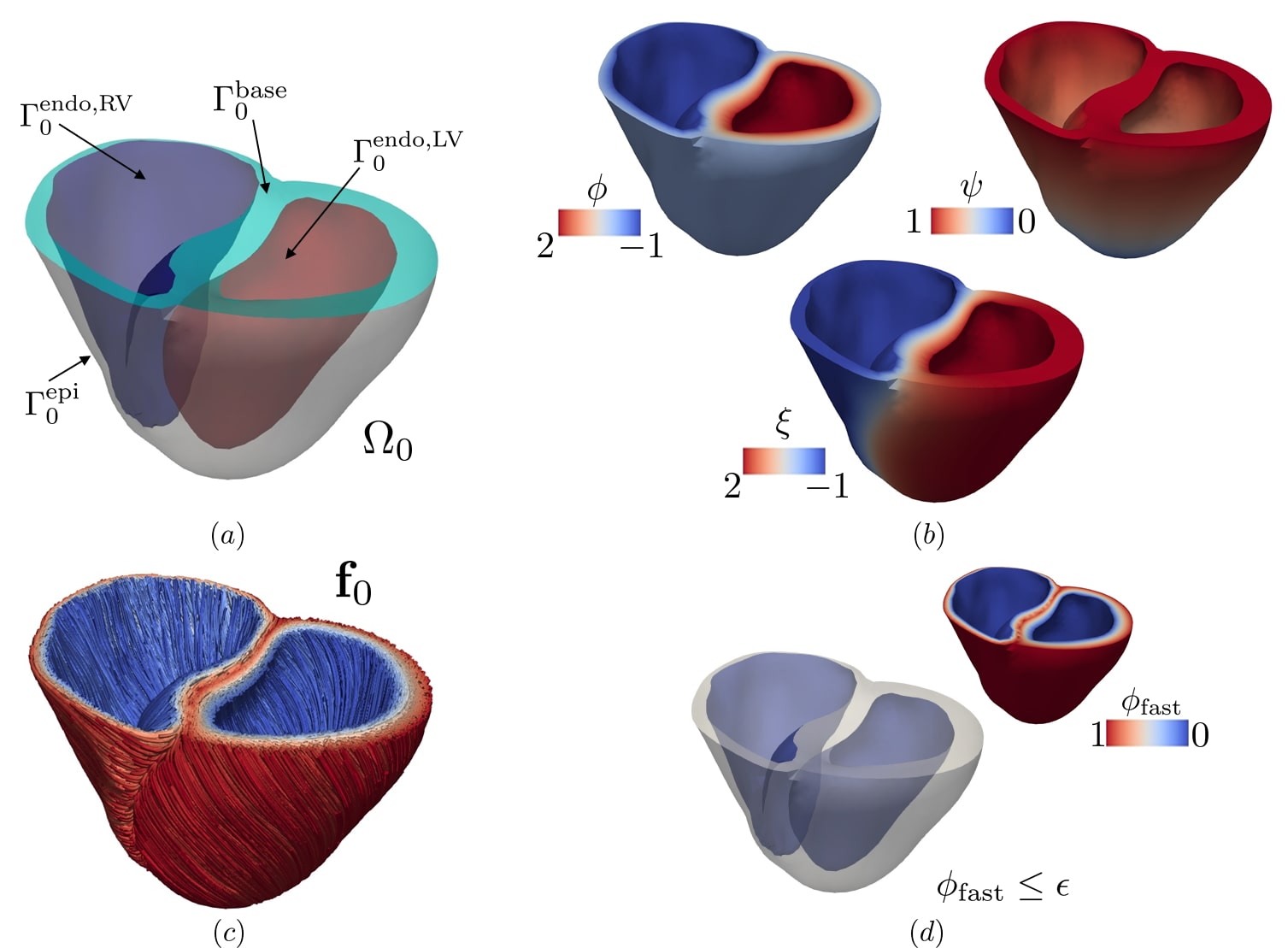}
	\caption{Top left (a): representation of a realistic biventricular computational domain $\Omega_0$ whose border is partitioned in $\GammaEpi$, $\GammaBase$, $\GammaEndoLV$ and $\GammaEndoRV$. Top Right (b): solutions of the Laplace problem \eqref{eqn: laplace} defining $\phi$ the transmural, $\psi$ the apico-basal and $\xi$ the inter-ventricular distances that are used to prescribe the myofiber orientations using LDRBM of type D-RBM. Bottom left (c): fiber field $\fZero$ obtained using D-RBM. Bottom  right (d): $\phi_{\mathrm{fast}}$ solution of the Laplace problem \eqref{eqn: laplace} used to build the fast endocardial layer $\phi_{\mathrm{fast}} \le \epsilon$~\cite{lee2019rule}.}
	\label{fig: tags_fibers_fast}
\end{figure}

We denote by $\Omega_0$ the computational domain in the reference configuration, see Figure~\ref{fig: tags_fibers_fast}(a), representing the region occupied by the left and right ventricles, whose boundary $\partial \Omega_0$ is partitioned into the epicardium $\GammaEpi$, the left $\GammaEndoLV$
and right $\GammaEndoRV$ endocardial surfaces and the biventricular base $\GammaBase$ (namely an artificial basal plane located well below the cardiac valves), so that we have $\partial \Omega_0=\GammaClosedEpi~\cup ~\GammaClosedEndoLV~\cup~\GammaClosedEndoRV~\cup~\GammaClosedBase$.
\subsection{Fibers generation}
\label{subsec: fibers}

To prescribe the cardiac muscle fiber architecture in the biventricular computational domain $\Omega_0$,  we use a particular class of Rule-Based-Methods (RBMs), known as Laplace–Dirichlet-Rule-Based-Methods (LDRBMs) \cite{bayer2012novel,wong2014generating,bayer2005laplace}. Specifically, we consider three LDRBMs, respectively proposed by Rossi et al. (R-RBM) \cite{quarteroni2017integrated}, Bayer et al. (B-RBM) \cite{bayer2012novel} and Doste et al. (D-RBM) \cite{doste2019rule}, that were recently reviewed in a communal mathematical description and extended to embed specific fiber directions for the right ventricle (RV) in \cite{piersanti2021modeling}.

LDRBMs define the transmural $\phi$ (from epicardium to endocardium), the apico-basal $\psi$ (from apex to basal plane) and the inter-ventricular $\xi$ (from the left to right endocardia) distances as the solutions of suitable Laplace boundary-value problems of the type
\begin{equation}
\label{eqn: laplace}
\begin{cases}
-\Delta \chi=0 &\qquad{\text{in }}\Omega_0,
\\
\chi = \chiA &\qquad{\text{on }}\GammaA,
\\
\chi = \chiB &\qquad{\text{on }}\GammaB,
\\
\nabla \chi \cdot \mecNref=0 &\qquad{\text{on }}\GammaN,
\end{cases}
\end{equation}
where $\chi=\phi, \psi, \xi$ denotes a generic unknown, $\chiA,\,\chiB \in \mathbb{R}$ are suitable boundary data set on generic partitions of the boundary $\GammaA,\,\GammaB,\,\GammaN$, with $\GammaClosedA \cup \GammaClosedB \cup \GammaClosedN=\partial\Omega_0$ and $\mecNref$ is defined as the outer normal vector, see Figures~\ref{fig: tags_fibers_fast}(b). For each point of the biventricular domain,  the transmural and apico-basal distances are used to build an orthonormal local coordinate axial system $[\eL, \eN, \eT]$ owing to $\eT= \tfrac{\nabla \phi}{\left\lVert \nabla \phi \right\rVert}$, $\eN=\tfrac{\nabla \psi - (\nabla \psi \cdot \eT )\eT}{\left\lVert \nabla \psi - (\nabla \psi \cdot \eT )\eT \right\rVert}$ and $\eL=\eN \times \eT$, defined as the unit transmural, longitudinal and normal directions, respectively. Finally, the reference frame $[\eL, \eN, \eT]$ is properly rotated with the purpose of defining the myofiber orientations:
\begin{equation*}
[\eL, \eN, \eT] \xrightarrow{\alphai, \betai} [\fZero, \nZero, \sZero],\quad \ii=\text{LV,RV},
\end{equation*}
where $\fZero$ is the fiber direction, $\nZero$ is the sheet-normal direction, $\sZero$ is the sheet direction,  $\mathrm{i}=\text{LV,RV}$ refers to LV or RV, and $\alphai$ and $\betai$ are suitable helical and sheetlet angles following linear relationships $\theta_{\ii}(d_{\ii}) = \thetaEpii(1-d_{\ii})+\thetaEndoi d_{\ii}$, (with $\theta=\alpha, \beta$ and $\mathrm{i}=\text{LV,RV}$) in which $d_{\ii}\in [0,1]$ is the transmural normalized distance and $\thetaEndoi$, $\thetaEpii$ are suitable prescribed rotation angles on the endocardium and epicardium, see Figure~\ref{fig: tags_fibers_fast}(c). To prescribe different myofiber orientations for LV and RV, we employ the inter-ventricular distance $\xi$ in which positive values of $\xi$ identify the LV, whereas negative values refer to the RV~\cite{piersanti2021modeling}. 
Moreover, we define the normalized inter-ventricular distance $\XiHat \in [0,1]$ by rescaling $\xi$. 

An example of LDRBM boundary-value solutions for the fiber generation procedure (of D-RBM type) is sketched in Figure~\ref{fig: tags_fibers_fast}(b). For further details about LDRBMs we refer to \cite{piersanti2021modeling}.

\subsection{3D-0D closed-loop electromechanical model}
\label{subsec: elettromech_model}

We provide a detailed description of the multiphysics and multiscale 3D biventricular EM model coupled with a 0D closed-loop (lumped parameters) hemodynamical model of the whole cardiovascular system, including the heart blood flow. Our model features several extensions and novel additions with respect to the work~\cite{regazzoni2020model,regazzoni2020numerical}, that is limited to the left ventricle. Our 3D-0D model is composed of four \textit{core models} supplemented by a suitable coupling condition between the 3D and the 0D model. The core models are related to the different interplaying physical processes (at the molecular, cellular, tissue and organ levels) involved in the heart pumping function: cardiac electrophysiology $(\PhyEP)$ \cite{levrero2020sensitivity,franzone2014mathematical,colli2018numerical,nobile2012active}, cardiomyocytes active contraction $(\PhyAct)$ \cite{rossi2014thermodynamically,niederer2011length,ruiz2014mathematical,land2017model,regazzoni2020biophysically,regazzoni2020machine},  tissue mechanics $(\PhyMec)$ \cite{guccione1991passive,guccione1991finite,ogden1997non,holzapfel2009constitutive} and blood circulation $(\PhyCirc)$ \cite{regazzoni2020model,hirschvogel2017monolithic,augustin2020physiologically,kerckhoffs2007coupling,blanco20103d,arts2005adaptation,quarteroni2016geometric}. The coupling condition is established by the volume conservation constraints~$(\PhyCoupl)$~\cite{regazzoni2020model}.

The model unknowns are:
\begin{equation*} 
\begin{aligned}
	u \colon &\Omega_0 \times (0,T)] \to \mathbb{R} , &
	\boldsymbol{w} \colon &\Omega_0 \times (0,T] \to \mathbb{R}^{\EPnumGating} , \\
	\ActStateHF \colon &\Omega_0 \times (0,T] \to \mathbb{R}^{\ActNumVariables} , &
	\displ \colon &\Omega_0 \times (0,T] \to \mathbb{R}^3, &
	\Circ \colon &(0,T] \to \mathbb{R}^{\CircNumVariables} , \\
	\PLV \colon &(0,T] \to \mathbb{R}, & \PRV \colon &(0,T] \to \mathbb{R} ,
\end{aligned}
\end{equation*} 
where $u$ is the transmembrane action potential, $\boldsymbol{w}$ the ionic variables vector, $\ActStateHF$ the state variables of the active force generation model, $\displ$ the tissue mechanical displacement, $\Circ$ the state vector of the circulation model (including pressures, volumes and fluxes of the different compartments composing the vascular network) and $\PLV$ and $\PRV$ are the left and right ventricular pressures, respectively.

Given the computational domain $\Omega_0$ and the time interval $t \in (0,T]$, our complete 3D-0D model reads:
\input{EMmodel_biv_CMAME.tex} 

The definition of the vectors $\BCmecVbaseLV$ and $\BCmecVbaseRV$, entering in the boundary conditions of the mechanical model $(\PhyMec)$ will be provided later.
Finally, the model is closed by the initial conditions in $\Omega_0 \times \{0\}$:
\begin{equation*}
u = u_0, \quad \boldsymbol{w} = \boldsymbol{w}_0, \quad \ActStateHF = \ActStateHF_0, \quad \displ = \displ_0, \quad \dfrac{\partial \displ}{\partial t} = \Dot{\displ}_0, \quad \Circ=\CircInit.
\end{equation*}

\subsubsection{Electrophysiology $(\PhyEP)$}
\label{subsec: electrophysiology}
We model the electric activity in the cardiac tissue by means of problem $(\PhyEP)$, that is the monodomain equation~\eqref{eqn: EP1} endowed with a suitable ionic model~\eqref{eqn: EP2} for the human ventricular action potential~\cite{levrero2020sensitivity,franzone2014mathematical,nobile2012active}. In the electrophysiology core model $(\PhyEP)$, the unknowns are the transmembrane potential $u$ and the ionic variables $\boldsymbol{w}$.
The vector $\boldsymbol{w}=\{w_1,w_2,...,w_{\EPnumGating}\}$ encodes the gating-variables (representing the fraction of open channels per unit area across the cell membrane) and the concentration of specific ionic species (among them the intracellular calcium ions concentration $\Cai$ plays a key role in the active force generation mechanism). The constant $\EPchim$ represents the surface area-to-volume ratio of cardiomyocytes, $\EPCm$ represents the trans-membrane capacitance per unit area. The applied current $\EPIapp$ mimics the effect of the Purkinje network \cite{vergara2014patient,costabal2016generating,landajuela2018numerical} modelled in this work by means of a surrogate fast endocardial conduction layer \cite{lee2019rule} represented by $\phi_{\mathrm{fast}}=\phi_{\mathrm{fast}}(\phi)\le \epsilon$ built as a function of the transmural distance defined in Section~\ref{subsec: fibers}, see also Figure~\ref{fig: tags_fibers_fast}(d).
The reaction terms $\EPIion$ and $\EPrhsGating$ (specified by the ionic model at hand) couple together the action potential propagation and the cellular dynamics.
Specifically, we use the ventricular ten Tusscher-Panfilov ionic model (TTP06, $\EPnumGating=18$), which is able to accurately describe ions dynamics across the cell membrane~\cite{ten2006alternans}. Furthermore, problem $(\PhyEP)$ is equipped with homogeneous Neumann boundary conditions~\eqref{eqn: EP3}.

The action potential propagation is driven by the diffusion term $\nabla \cdot ( J \mecF^{-1} \EPdiffTens \mecF^{-T} \nabla u)$ where we introduced the deformation gradient tensor $\mecF = \identity + \nabla \displ$ with $J = \det(\mecF)>0$. The diffusion tensor reads:
\begin{equation*} 
\EPdiffTens = \sigma_{\ell} (\phi_{\mathrm{fast}}) \frac{\mecF\fZero \otimes \mecF\fZero}{\|\mecF \fZero \|^2} + \sigma_{\text{t}} (\phi_{\mathrm{fast}}) \frac{\mecF\sZero \otimes \mecF\sZero}{\|\mecF \sZero \|^2} + \sigma_{\text{n}} (\phi_{\mathrm{fast}}) \frac{\mecF\nZero \otimes \mecF\nZero}{\|\mecF \nZero \|^2},
\end{equation*}
where $\sigma_{\ell}(\phi_{\mathrm{fast}}), \sigma_{\text{t}}(\phi_{\mathrm{fast}}), \sigma_{\text{n}}(\phi_{\mathrm{fast}})$ are the longitudinal, transversal and normal conductivities, respectively, defined as
\begin{equation*}
	\sigma_{\text{k}}(\phi_{\mathrm{fast}})=
	\begin{cases}
		\sigma_{\text{k,fast}}
		& $if$ \; \phi_{\mathrm{fast}} \le \epsilon, \quad \mathrm{k=\ell,t,n},  \\
		\sigma_{\text{k,myo}}
		& $if$ \; \phi_{\mathrm{fast}} > \epsilon,  \quad \mathrm{k=\ell,t,n},
	\end{cases}
\end{equation*}
where $\sigma_{\text{k,fast}}$ and $\sigma_{\text{k,myo}}$ (with $\mathrm{k=\ell,t,n}$) are the prescribed conductivities inside and outside the fast endocardial layer, respectively.
\subsubsection{Activation $(\PhyAct)$}
\label{subsec: mechanicalactivation}
Mechanical activation of cardiac tissue is modeled by means of equation~\eqref{eqn: Act}, a system of ODEs standing for an Artificial Neural Network (ANN) based model that surrogates the so called RDQ18 high-fidelity model proposed in \cite{regazzoni2020machine}. The RDQ18 model is based on a biophysically detailed description of the microscopic active force generation mechanisms taking place at the scale of sarcomeres \cite{bers2001excitation,regazzoni2020biophysically}. The RDQ18-ANN model has the great advantage of strikingly reducing the computational burden associated to the numerical solution of the RDQ18 model, yet reproducing its results with a very good accuracy \cite{regazzoni2020machine}.

In the activation core model $(\PhyAct)$ the unknown is the two-variable state vector $\ActStateHF$. The input variables
are the intracellular calcium ions concentration $\Cai$, provided by the TTP06 ionic model, and the sarcomere length $\SL$ defined as $\SL=\SL_0 \sqrt{\IIVf(\displ)}$, 
where $\SL_0$ denotes the sarcomere length at rest and $\IIVf = \mecF \fZero \cdot \mecF \fZero$ is a measure of the tissue stretch along the fibers direction. This creates a feedback between the mechanical model $(\PhyMec)$ and the force generation model $(\PhyAct)$ \cite{regazzoni2020model}.

The RDQ18-ANN output is the permissivity $P \in [0,1]$ which is obtained as a function of $\ActStateHF$: $P=G(\ActStateHF)$ where $G$ is a linear function defined in \cite{regazzoni2020machine}. Since $P$ is the fraction of the contractile units in the force-generation state, the active tension is given by  $\Tens = \TensMax \, P$, where $\TensMax$ denotes the tension generated when all the contractile units are generating force (i.e. for $P=1$). Finally, to account for a different active tension between LV and RV we define a spatial heterogeneous active tension
\begin{equation*}
\Tens(\ActStateHF) = \TensMax G(\ActStateHF) \,\left[ \XiHat + \Cr ( 1 - \XiHat ) \right],
\end{equation*}
where $\XiHat \in [0,1]$ is the normalized inter-ventricular distance, defined in Section~\ref{subsec: fibers}, and $\Cr\in (0,1]$ represents the left-right ventricle contractility ratio.

\subsubsection{Mechanics $(\PhyMec)$}
\label{subsec: actpassmechanics}
The mechanical response of the cardiac tissue is described by problem $(\PhyMec)$ under the hyperelasticity assumption and by adopting an active stress approach \cite{guccione1991finite,ogden1997non}. The unknown is the displacement $\displ$, whereas $\rho_{\text{s}}$ is the density. The first Piola-Kirchhoff stress tensor $\tenspiola = \tenspiola(\displ, \Tens)$ is additively decomposed according to
\begin{equation}\label{eqn: piola}
\tenspiola(\displ, \Tens)
= \dfrac{\partial \mathcal{W}(\mecF)}{\partial \mecF}
+ \Tens(\XiHat,\ActStateHF) \bigg[ n_{\text{f}} \frac{\mecF \fZero \otimes \fZero}{ \sqrt{\IIVf}} + n_{\text{s}}\frac{\mecF \sZero \otimes \sZero}{ \sqrt{\IIVs}} + n_{\text{n}}\frac{\mecF \nZero \otimes \nZero}{ \sqrt{\IIVn}}\bigg],
\end{equation}
where the first term represents the passive mechanics with $\mathcal{W}: \LinPlus \to \mathbb{R}$ being the strain energy density function, whereas the second one stands for the orthotropic active stress, with $\Tens(\XiHat,\ActStateHF)$ the active tension provided by the activation model $(\PhyAct)$. Moreover, $\IIVs = \mecF \sZero \cdot \mecF \sZero$ and $\IIVn = \mecF \nZero \cdot \mecF \nZero$ are the tissue stretches along the sheet and sheet-normal directions, respectively, and $n_{\text{f}}$, $n_{\text{s}}$ and $n_{\text{n}}$ the prescribed proportion of active tension along the fiber, sheet and sheet-normals directions, respectively. Notice that the orthotropic active stress tensor~\eqref{eqn: piola} surrogates the contraction caused by dispersed myofibers~\cite{guan2020effect,guan2021modelling,genet2014distribution,sack2018construction}.

To model the passive behaviour of the cardiac tissue, we employ the orthotropic Guccione constitutive law \cite{guccione1991passive}, according to which the strain energy function is defined as
\begin{equation*}
\mathcal{W} = \dfrac{\kappa}{2} \left( J - 1 \right) \text{log}(J) + \dfrac{a}{2} \left( e^Q  - 1 \right),
\end{equation*}
where the first term is the volumetric energy with the bulk modulus $\kappa$, which penalizes large variation of volume to enforce a weakly incompressible behaviour \cite{peng1997compressible,doll2000development}, and the exponent $Q$ reads
\begin{equation*}
\begin{split}
Q &= b_{\text{ff}} E_{\text{ff}}^2  + b_{\text{ss}} E_{\text{ss}}^2 + b_{\text{nn}} E_{\text{nn}}^2+ b_{\text{fs}} \left( E_{\text{fs}}^2 + E_{\text{sf}}^2 \right) + b_{\text{fn}} \left( E_{\text{fn}}^2 + E_{\text{nf}}^2 \right) + b_{\text{sn}} \left( E_{\text{sn}}^2 + E_{\text{ns}}^2 \right),
\end{split}
\end{equation*}
where $a$ is the stiffness scaling parameter, $E_\text{ij} = \boldsymbol{E} \boldsymbol{\mathrm{i}}_\text{0} \cdot \boldsymbol{\mathrm{j}}_\text{0}$, for  $\mathrm{i, j} \in \{ \mathrm{f, s, n} \}$ and $\boldsymbol{\mathrm{i}}_\text{0},\boldsymbol{\mathrm{j}}_\text{0} \in \{ \mathrm{\fZero, \sZero, \nZero} \}$,
are the entries of $\textbf{E} = \tfrac{1}{2} \left( \mecC - \identity \right)$, i.e the Green-Lagrange strain tensor, being $\mecC = \mecF^{T} \mecF$ the right Cauchy-Green deformation tensor.

To model the mechanical constraint provided by the pericardium \cite{Gerbi2018monolithic,pfaller2019importance,strocchi2020simulating}, we impose at the epicardial boundary $\GammaEpi$ a generalized Robin boundary condition~\eqref{eqn: Mec2} by defining the tensors $		\BCmecKepiTens =
\BCmecKepiT (\mecNref \otimes \mecNref - \identity) -
\BCmecKepiN (\mecNref \otimes \mecNref)$ and $		\BCmecCepiTens =
\BCmecCepiT (\mecNref \otimes \mecNref - \identity) -
\BCmecCepiN (\mecNref \otimes \mecNref),$
where $\BCmecKepiN$, $\BCmecKepiT$, $\BCmecCepiN$, $\BCmecCepiT \in \mathbb{R}^+$ are the stiffness and viscosity parameters of the epicardial tissue in the normal and tangential directions, respectively. Normal stress boundary conditions ~\eqref{eqn: Mec3}--\eqref{eqn: Mec4} were imposed at the endocardia $\GammaEndoLV$ and $\GammaEndoRV$ of both ventricles where $\PLV(t)$ and $\PRV(t)$ represent the pressure exerted by the blood in the left and right ventricular chambers, respectively.
To take into account the effect of the neglected part, over the basal plane, on the biventricular domain, we set on $\GammaBase$ the energy-consistent boundary condition~\eqref{eqn: Mec5} in weighted-stress-distribution form, where
\begin{equation}\label{eqn: weighted_bc}
\BCmecVbaseLV(t,\XiHat) = \XiHat \displaystyle \frac{\int_{\GammaEndoLV} \, J \mecF^{-T} \mecNref d \Gamma_0}{\int_{\GammaBase} \XiHat \, | J \mecF^{-T} \mecNref | d \Gamma_0},
\qquad 
\BCmecVbaseRV(t,\XiHat) = (1-\XiHat) \displaystyle \frac{\int_{\GammaEndoRV} \, J \mecF^{-T} \mecNref d \Gamma_0}{\int_{\GammaBase} (1-\XiHat) \, | J \mecF^{-T} \mecNref | d \Gamma_0}.
\end{equation}
The energy-consistent boundary condition considered in this work is the extension to the biventricular case of the energy-consistent boundary condition originally proposed in \cite{regazzoni2020machine} for LV.
The complete derivation can be found in Appendix~\ref{app:BC}.

\subsubsection{Blood circulation $(\PhyCirc)$ and coupling conditions $(\PhyCoupl)$}
\begin{figure}
	\centering
	\includegraphics[keepaspectratio, width=0.5\textwidth]{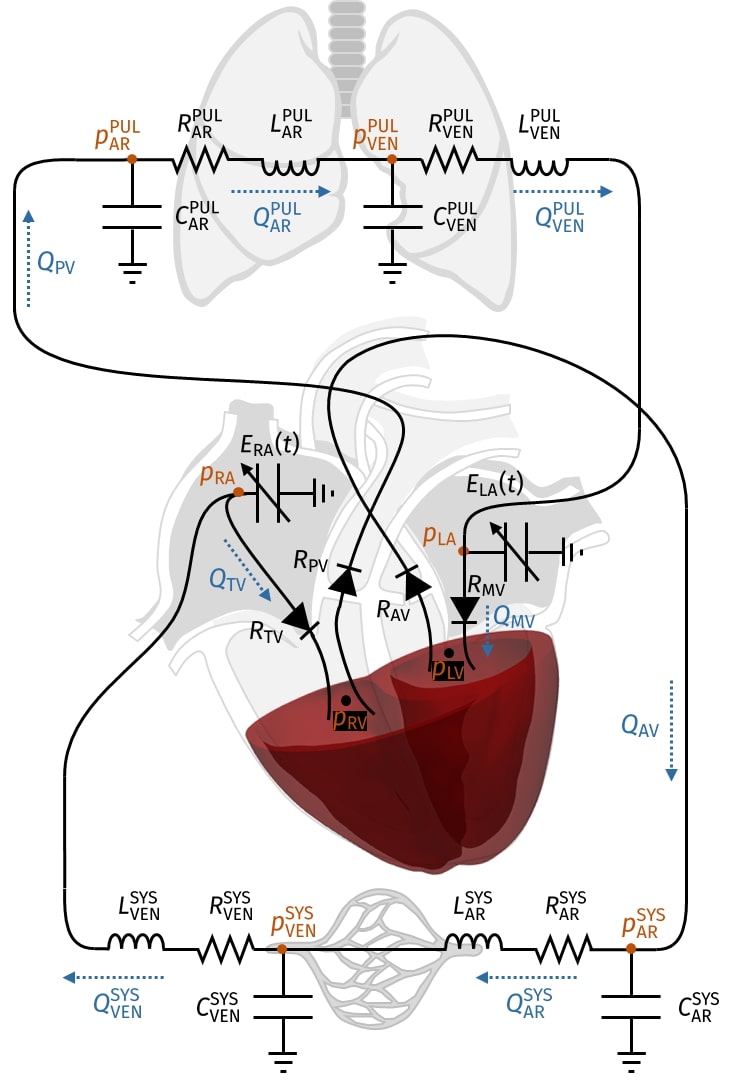}
	\caption{3D-0D coupling between the biventricular 3D EM model and the 0D circulation model.}
	\label{fig: circulation}
\end{figure}
\label{sec: mathematicalmodeling: circulation}
We model the blood circulation through the entire cardiovascular system (i.e. equation~\eqref{eqn: Circ}) by means of a closed-loop model, inspired by~\cite{hirschvogel2017monolithic,blanco20103d} and recently proposed in \cite{regazzoni2020model}. In the 0D closed-loop model, systemic and pulmonary circulations are modeled with RLC circuits, heart chambers are described by time-varying elastance elements and non-ideal diodes stand for the heart valves \cite{regazzoni2020model}.

The circulation core model $(\PhyCirc)$ is represented by a system of ODEs expressed by equation~\eqref{eqn: Circ}, where $\CircRhs$ is a proper function (defined in \cite{regazzoni2020model}) and $\Circ(t)$ includes pressures, volumes and fluxes of the different compartments composing the vascular network:
\begin{equation*}
	\begin{split}
		\Circ(t) =
		(&\VLA(t), \VLV(t), \VRA(t), \VRV(t), \ParSYS(t), \PvnSYS(t), \ParPUL(t), \PvnPUL(t), \\&\QarSYS(t), \QvnSYS(t), \QarPUL(t), \QvnPUL(t))^T.
	\end{split}
\end{equation*}
Here $\VLA$, $\VRA$, $\VLV$ and $\VRV$ refer to the volumes of left atrium, right atrium, LV and RV, respectively;
$\ParSYS$, $\QarSYS$, $\PvnSYS$, $\QvnSYS$, $\ParPUL$, $\QarPUL$, $\PvnPUL$ and $\QvnPUL$ express pressures and flow rates of the systemic and pulmonary circulation (arterial and venous). For the complete mathematical description of the 0D circulation lumped model we refer to \cite{regazzoni2020model}.
To couple the 0D circulation model $(\PhyCirc)$ with the 3D biventricular model, given by $(\PhyEP)$--$(\PhyAct)$--$(\PhyMec)$, we follow the strategy proposed in \cite{regazzoni2020model}:
we replace the time-varying elastance elements representing LV and RV in the circulation model with its corresponding 3D electromechanical description, obtaining the coupled 3D-0D model depicted in Figure~\ref{fig: circulation}. With this aim, we introduce the volume-consistency coupling conditions $(\PhyCoupl)$ where
\begin{equation*}
	\displaystyle
	\Vthreedimi(\displ(t)) =
	\int_{\GammaEndoi} J(t)
	\left(\left( \mathbf{h} \otimes \mathbf{h}  \right) \left(\mathbf{x} + \displ(t) - \mathbf{b}_i \right) \right)
	\cdot \mecF^{-T}(t) \mecNref \, d \Gamma_0 \qquad \mathrm{i=LV,RV}
\end{equation*}
wherein $\mathbf{h}$ is a vector orthogonal to LV/RV centreline (i.e. lying on the biventricular base) and $\mathbf{b}_i$ lays inside LV/RV \cite{regazzoni2020model}.

In virtue of the introduced conditions $(\PhyCoupl)$, in the 3D-0D coupled model~\eqref{eqn: EM} the left $\PLV(t)$ and right $\PRV(t)$ ventricular pressures are not determined by the 0D circulation model~equation~\eqref{eqn: Circ}, but rather act as Lagrange multipliers associated to the constraint $(\PhyCoupl)$.
\subsection{Reference configuration and initial tissue displacement}
\label{subsec: refconf}
Cardiac geometries are acquired from in vivo medical images through imaging techniques. These geometries are in principle not stress free, due to the blood pressure acting on the endocardia. Therefore, we need to estimate the unloaded (i.e. stress-free) configuration (also named \textit{reference configuration}) to which the 3D-0D model~\eqref{eqn: EM} refers. To recover the reference configuration $\Omega_0$, starting from a geometry $\OmegaRef$ acquired from medical images (typically during the diastolic phase), we extend to the biventricular case the procedure proposed for LV in \cite{regazzoni2020model}.

We assume that the configuration $\OmegaRef$ is acquired during the diastole, when the biventricular geometry  is loaded with $\PLV = \pressLVRef$, $\PRV = \pressRVRef$ and a residual active tension $\Tens = \tensRef>0$ is present. To recover the reference configuration $\Omega_0$ we solve the following inverse problem: find the domain $\Omega_0$ such that, if we inflate $\Omega_0$ by $\displ$, solution of the differential problem\footnote{The problem~\eqref{eqn: mechanics_steadystate} is derived from $(\PhyMec)$ setting aside the time dependent terms.}
\begin{equation}
\begin{cases}\label{eqn: mechanics_steadystate}
	\nabla \cdot \tenspiola(\displ, \Tens) = \boldsymbol{0} & $in$ \; \Omega_0 , \\
	\tenspiola(\displ, \Tens) \mecNref +
	\BCmecKepiTens \displ
	= \mathbf{0}
	& $on$ \; \GammaEpi , \\
\tenspiola(\displ, \Tens) \mecNref = -\PLV(t) \, J \mecF^{-T} \mecNref & $on$ \; \GammaEndoLV, \\
\tenspiola(\displ, \Tens) \mecNref = -\PRV(t) \, J \mecF^{-T} \mecNref & $on$ \; \GammaEndoRV, \\
\tenspiola(\displ, \Tens) \mecNref = \displaystyle | J \mecF^{-T} \mecNref | \left[ \PLV \BCmecVbaseLV(\XiHat) + \PRV \BCmecVbaseRV(\XiHat) \right]
& $on$ \; \GammaBase,
\end{cases}
\end{equation}
obtained for $\PLV = \pressLVRef$, $\PRV = \pressRVRef$ and $\Tens = \tensRef$, we get the domain $\OmegaRef$.

After recovering $\Omega_0$, we inflate the biventricular reference configuration $\Omega_0$ by solving again problem~\eqref{eqn: mechanics_steadystate}, where we set the pressures $\PLV=p_\text{LV}^\text{ED}$ and $\PRV=p_\text{RV}^\text{ED}$ with the superscript ED stands for the end-diastolic phase. The values $p_\text{LV}^\text{ED}$ and $p_\text{RV}^\text{ED}$ are chosen to bring the biventricular domain to defined end diastolic volumes for the left $V_\text{LV}^\text{ED}$ and right $V_\text{RV}^\text{ED}$ ventricles. In this way we obtain the end-diastolic configuration for the biventricular geometry. Hence, the solution $\displ$ of the problem~\eqref{eqn: mechanics_steadystate} is set as initial condition $\displ_0$ for $\displ$ in $(\PhyMec)$. The above procedure is represented in step 4 of Figure~\ref{fig: pipeline}.

%% file: parts_numerics.tex
\section{Numerical approximation}
\label{sec: numericaldiscretization}
In this section we illustrate the numerical discretization of the different core models composing the 3D-0D problem~\eqref{eqn: EM} along with the strategy that we adopt to reach a limit-cycle.

\subsection{Space and time discretizations}
\label{sec:space_time_discretization}
\begin{figure}[t!]
	\centering
	\includegraphics[keepaspectratio, width=1\textwidth]{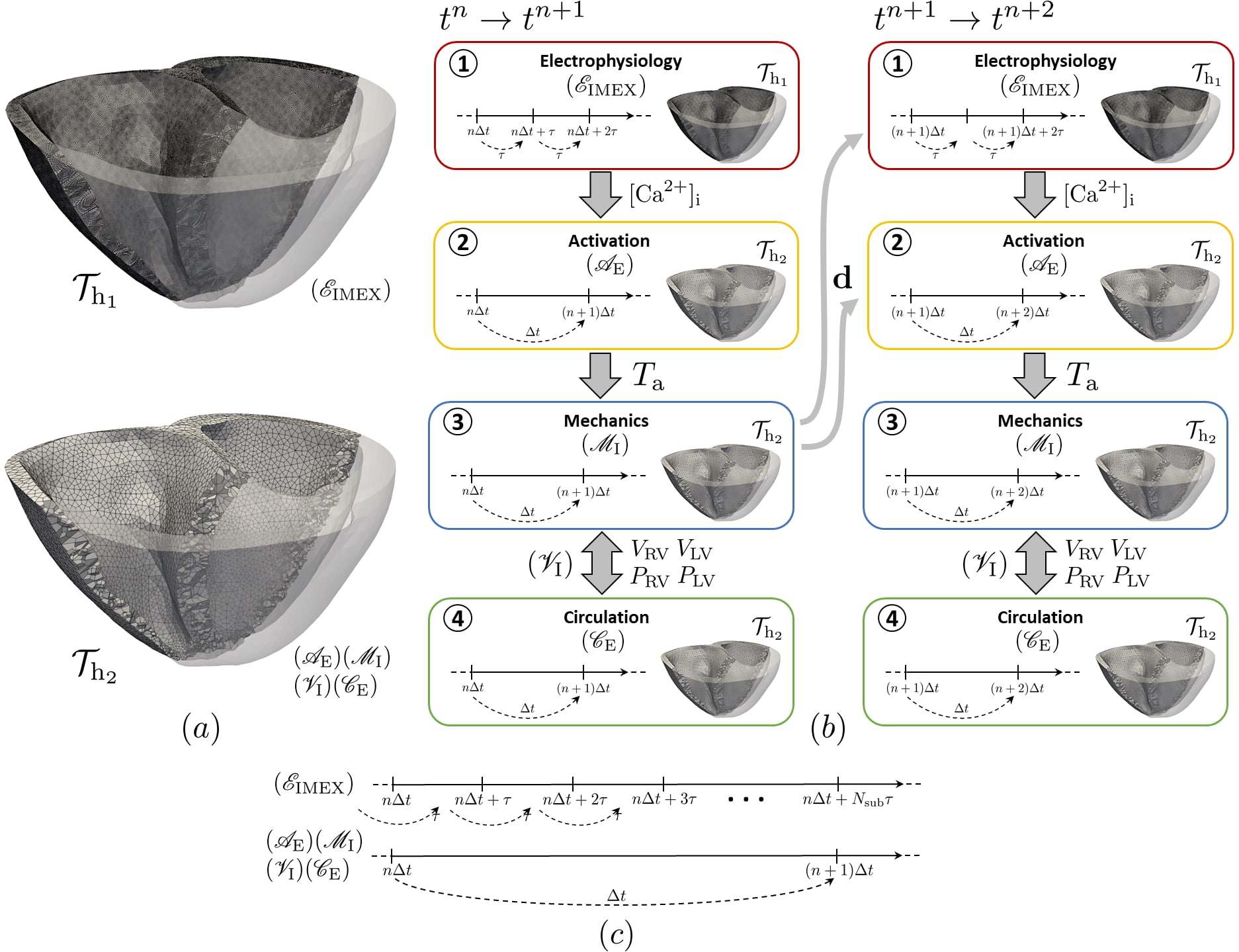}
	\caption{Segregated-intergrid-staggered numerical scheme: (a) nested meshes $\mathcal{T}_{\hone}$ and $\mathcal{T}_{\htwo}$ (with $\hone<\htwo$); (b) schematic representation of the numerical scheme; (c) graphical representation of the time advancement scheme.}
	\label{fig: numerical_scheme}
\end{figure}

For the numerical approximation of the 3D-0D coupled model~\eqref{eqn: EM}
we follow the approach proposed in~\cite{regazzoni2020numerical}, which is extended here to the biventricular case. The core models $(\PhyEP)-(\PhyAct)-(\PhyMec)-(\PhyCirc)$ are solved sequentially in a segregated manner by using different resolutions in space and time, to properly handle the different space and time scales of the core models contributing to both cardiac EM and blood circulation~\cite{nordsletten2011coupling,quarteroni2017integrated,quarteroni2019cardiovascular}. For this reason we call this numerical approach Segregated-Intergrid-Staggered (SIS). 

For the space discretization, we use the Finite Element Method (FEM) with continuous Finite Elements (FEs) of order 1 ($\mathbb{Q}_1$) and hexahedral meshes \cite{quarteroni2009numerical}. We consider two nested meshes $\mathcal{T}_{\hone}$ and $\mathcal{T}_{\htwo}$ of the computational domain $\Omega_0$ ($\hone$ and $\htwo$, with $\hone<\htwo$, represent the mesh sizes), where $\mathcal{T}_{\hone}$ is built by uniformly refining $\mathcal{T}_{\htwo}$~\cite{Africa2019,BursteddeWilcoxGhattas11}, see Figure~\ref{fig: numerical_scheme}(a). We adopt the finer mesh $\mathcal{T}_{\hone}$ for $(\PhyEP)$, where it is essential to accurately capture the dynamics of travelling waves, while the coarser one ($\mathcal{T}_{\htwo}$) is used for both $(\PhyAct)$ and $(\PhyMec)$ \cite{augustin2016anatomically,regazzoni2020numerical,colli2018numerical}.
We employ an efficient intergrid transfer operator between the nested grids $\mathcal{T}_{\hone}$ and $\mathcal{T}_{\htwo}$, which allows to evaluate the feedback between $(\PhyEP)$ and $(\PhyAct)-(\PhyMec)$~\cite{regazzoni2020numerical}.
In \cite{regazzoni2020numerical}, the displacement field $\displ$ is interpolated on $\mathcal{T}_{\hone}$ and $\nabla \displ$ is assembled on the fine mesh directly.
Here, we follow the more effective strategy proposed in \cite{salvador2020intergrid}, where $\nabla \displ$ is recovered on $\mathcal{T}_{\htwo}$ thanks to an $L^2$ projection \cite{Africa2019}.
Then, $\nabla \displ$ is interpolated on $\mathcal{T}_{\hone}$.

For the time discretization, we use Finite Difference schemes~\cite{quarteroni2010numerical}.
The cardiac electrophysiology model is solved by means of the Backward Differentiation Formula of order 2 (BDF2).
We adopt an implicit-explicit (IMEX) scheme, denoted by $(\PhyEPIMEX)$, where the diffusion term is treated implicitly, whereas the ionic and reaction terms explicitly~\cite{regazzoni2020numerical,niederer2011simulating}.
For both mechanical activation and passive mechanics we employ the BDF1 scheme, where $(\PhyActE)$ is advanced in time with an explicit method, whereas a fully implicit scheme is used for $(\PhyMecI)-(\PhyCouplI)$ \cite{regazzoni2020numerical}.
Finally, we employ an explicit $4^{th}$ order Runge-Kutta method (RK4) for $(\PhyCircE)$ \cite{regazzoni2020numerical}.

We use two different time steps, $\Delta t$ for $(\PhyActE)-(\PhyMecI)-(\PhyCouplI)-(\PhyCircE)$ and $\tau = \Delta t/\Nsub$ for $(\PhyEPIMEX)$, with $\Nsub \in \mathbb{N}$, see Figure~\ref{fig: numerical_scheme}(c).
We first update the variables of $(\PhyEPIMEX)$, then those of $(\PhyActE)$ and finally, after updating the unknowns of $(\PhyMecI)-(\PhyCouplI)$, we update the ones of $(\PhyCircE)$, see Figure~\ref{fig: numerical_scheme}(b).

The whole algorithm for the SIS numerical scheme is reported in Figure~\ref{fig: numerical_scheme}.

\subsection{3D-0D coupled problem resolution}
\label{sec: saddle-point}

We couple the 3D mechanical model $(\PhyMec)$ with the 0D closed-loop hemodynamical model $(\PhyCirc)$ by means of the volume conservation constraints~$(\PhyCoupl)$, where the pressures of LV and RV act as Lagrange multipliers~\cite{regazzoni2020numerical}.
In Figure~\ref{fig: numerical_scheme}(b) (steps 3-4) we obtain a saddle point problem $(\PhyMecI)-(\PhyCouplI)$.

We introduce the discrete times $t^n=n\Delta t$, $n\ge0$ and we denote by $\adofs_{\text{h}}^{n} \simeq \adofs_{\text{h}}(t^n)$ the fully discretized FEM approximation of the generic (scalar, vectorial or tensorial) variable $\adofs(t)$ (i.e. the vector collecting the DOFs defined over the computational mesh $\mathcal{T}_{\htwo}$ at time $t^n$).
Then, at each time step $t^{n+1}$, the fully discretized version of $(\PhyMecI)-(\PhyCouplI)$ reads:
\begin{equation}
\begin{cases}\label{eqn: Msi:biventricle}
&\left( \rho_s \dfrac{1}{(\Delta t)^2} \mathcal{M} + \dfrac{1}{\Delta t} \mathcal{F} + \mathcal{G} \right) \displdofs_{\text{h}}^{n+1} + \boldsymbol{S}(\displdofs_{\text{h}}^{n+1}, \TensSubscriptDofs{\text{h}}^{n+1}) \\
&\qquad\qquad= \rho_s \dfrac{2}{(\Delta t)^2} \mathcal{M} \displdofs_{\text{h}}^n - \rho_s \dfrac{1}{(\Delta t)^2} \mathcal{M} \displdofs_{\text{h}}^{n-1} + \dfrac{1}{\Delta t} \mathcal{F} \displdofs_{\text{h}}^n
\\&\qquad\qquad
+ \PLV^{n+1}\boldsymbol{P}_{\text{LV}}(\displdofs_{\text{h}}^n,\displdofs_{\text{h}}^{n+1})
+ \PRV^{n+1}\boldsymbol{P}_{\text{RV}}(\displdofs_{\text{h}}^n,\displdofs_{\text{h}}^{n+1}) \\
& \VLV(\Circ^{n+1}) = \VLVthreedim(\displdofs_{\text{h}}^{n+1}) \\
& \VRV(\Circ^{n+1}) = \VRVthreedim(\displdofs_{\text{h}}^{n+1})
\end{cases}
\end{equation}
where we introduced
\begin{equation*}
	\displaystyle \mathcal{M}_{\text{ij}} = \int_{\Omega_0} \boldsymbol{\phi}_\text{j} \cdot \boldsymbol{\phi}_\text{i} \, d \Omega_0,
	\;\;\;\;
	\boldsymbol{S}_{\text{i}}=\int_{\Omega_0} \tenspiola(\displdofs_{\text{h}}^n, \TensSubscriptDofs{\text{h}}^n):\nabla \boldsymbol{\phi}_\text{i} \, d \Omega_0,
\end{equation*}
\begin{equation*}
	\mathcal{F}_{\text{ij}} = \int_{\GammaEpi} \left[ \BCmecCepiT(\mecNref_{\text{h}} \otimes \mecNref_{\text{h}} - \identity_{\text{h}}) - \BCmecCepiN (\mecNref_{\text{h}} \otimes \mecNref_{\text{h}})  \right] \boldsymbol{\phi}_\text{j} \cdot \boldsymbol{\phi}_\text{i} \, d \Gamma_0,
\end{equation*}

\begin{equation*}
	\mathcal{G}_{\text{ij}} = \int_{\GammaEpi} \left[ \BCmecKepiT(\mecNref_{\text{h}} \otimes \mecNref_{\text{h}} - \identity_{\text{h}}) - \BCmecKepiN (\mecNref_{\text{h}} \otimes \mecNref_{\text{h}})  \right] \boldsymbol{\phi}_\text{j} \cdot \boldsymbol{\phi}_\text{i} \, d \Gamma_0,
\end{equation*}

\begin{equation*}
	\begin{split}
		\boldsymbol{P}_\text{k,i} &=  \int_{\GammaBase}
		| J^{n+1}_{\text{h}} (\mecF^{n+1}_{\text{h}})^{-T} \mecNref_{\text{h}} |\BCmecVbaseDiscr{\text{k},\text{h}}\cdot \boldsymbol{\phi}_\text{i} \, d \Gamma_0 \\
		&- \int_{\GammaEndo} J^{n+1}_{\text{h}} (\mecF^{n+1}_{\text{h}})^{-T} \mecNref_{\text{h}} \cdot \boldsymbol{\phi}_\text{i} \, d \Gamma_0, \qquad \text{k=LV, RV.}
	\end{split}
\end{equation*}
Here $\mecF^{n+1}_{\text{h}}= \identity_{\text{h}} + \nabla \displdofs_{\text{h}}^{n+1}$ with $J^{n+1}_{\text{h}} = \det(\mecF^{n+1}_{\text{h}})$, $\{\boldsymbol{\phi}_\text{i}\}_{\text{i=1}}^{N_\text{d}}$ represents the set of basis functions for the finite dimensional space $[\mathcal{X}_{\text{h}}^s]^3$ with $\mathcal{X}_{\text{h}}^s= \{ v \in C^0(\Bar{\Omega}_0) : v|_\text{K} \in \mathbb{Q}_\text{s}(K), \, s \geq 1, \;\; \forall K \in \mathcal{T}_{\htwo} \}$, where $\mathbb{Q}_\text{s}(K)$ stands for the set of polynomials with degree smaller than or equal to $s$ over a mesh element $K$ and $N_{\displ} = \operatorname{dim}([\mathcal{X}_{\text{h}}^s]^3)$ is the numbers of DOFs for the displacement.

Moving all the terms to the right hand side, equation~\eqref{eqn: Msi:biventricle} can be compactly written as:
\begin{equation}\label{eqn: saddle-point}
\begin{cases}
\mathbf{r}_{\displ}(\displdofs_{\text{h}}^{n+1},\PLV^{n+1},\PRV^{n+1}) &= \mathbf{0} \\
r_{\PLV}(\displdofs_{\text{h}}^{n+1}) &= 0 \\
r_{\PRV}(\displdofs_{\text{h}}^{n+1}) &= 0 \\
\end{cases}
\end{equation}
for suitable functions $r_{\PLV}$, $r_{\PRV}$ and $\mathbf{r}_{\displ}$.

We solve the non-linear saddle-point problem~\eqref{eqn: saddle-point} by means of the Newton algorithm using, at the algebraic level, the Schur complement reduction \cite{benzi2005numerical}. More details about the resolution of the problem~\eqref{eqn: saddle-point} can be found in Appendix~\ref{app:schur}.

\subsection{Finding initial conditions for the multiphysics problem}
\label{subsec: initial_conditions}
\begin{figure}[t!]
	\centering
	\includegraphics[keepaspectratio, width=0.9\textwidth]{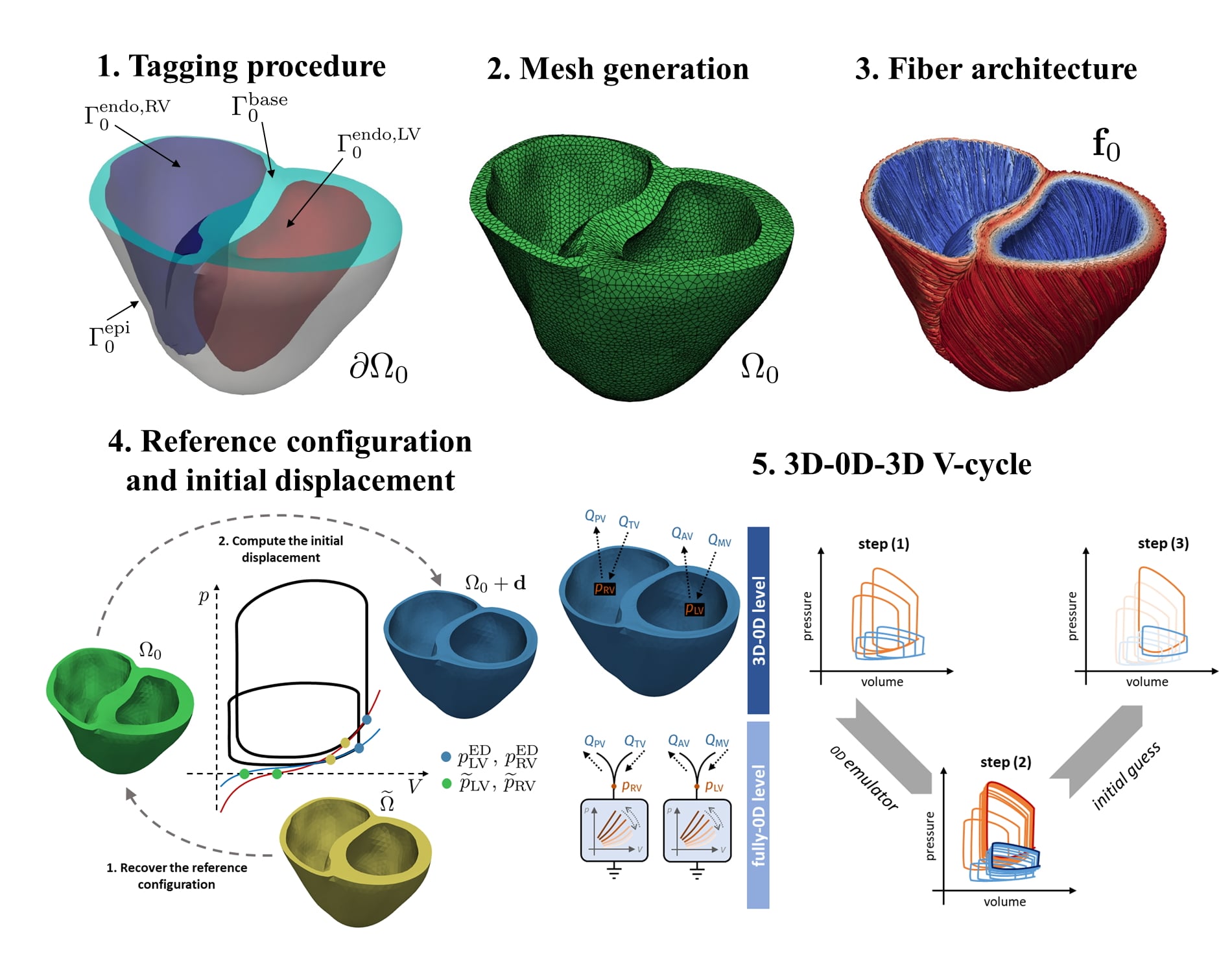}
	\caption{Graphical display of the whole pipeline for the initialization of a numerical simulation employing the 3D-0D EM model.}
	\label{fig: pipeline}
\end{figure}
The numerical results of the 3D-0D EM model typically feature a temporal transient, which lasts for several heartbeats and converges to a periodic solution, known as \textit{limit cycle}.
The outputs of clinical interest should be computed from the numerical solution that is associated with the limit cycle.
To reduce the computational overhead of reaching a periodic solution, we follow the strategy proposed in \cite{regazzoni2021emulator}, aimed at accelerating the convergence towards the limit cycle.
This strategy -- named 3D-0D-3D V-cycle -- comprises three stages (see point 5 of Figure~\ref{fig: pipeline}).
In a first step, three heartbeats are simulated with the 3D-0D model.
Then, based on the PV-loops obtained from the previous 3D-0D model, a 0D emulator of each ventricle is
built with the aim of surrogating the pressure-volume relationships, and substituted to the 3D model.
These emulators, coupled with the 0D model of blood circulation for the remaining compartments, allow to simulate the transient phase toward a periodic solution in less than one minute of computational time on a standard laptop.
Finally, the state obtained with this fully 0D model is used to initialize the 3D-0D model, and three additional heartbeats are simulated.
Overall, the computational cost of reaching the limit cycle amounts to that of simulating six heartbeats, regardless of the number of cycles required to converge to a periodic solution.
As a matter of fact, the computational time required by the 0D surrogate model is negligible compared to that of the full-order 3D-0D model.
More details on this pipeline are available in~\cite{regazzoni2021emulator}.

To find an initial guess for the remaining variables, we initialize the ionic model by running a 1000-cycle long single-cell simulation.
Similarly, we initialize the force generation model by means of a single-cell simulation with a constant calcium input (corresponding to the final calcium concentration of the single-cell ionic simulation) and a reference sarcomere length $SL = \SI{2.2}{\micro\meter}$.

The whole pipeline for the initialization of a numerical simulation employing the 3D-0D biventricular EM model is sketched in Figure~\ref{fig: pipeline}.

%% file: parts_results.tex
\section{Numerical results}
\label{sec: numericalresults}
In this section, we present several biventricular electromechanical simulations that employ the 3D-0D model discussed in Sections~\ref{sec: mathematicalmodeling} and~\ref{sec: numericaldiscretization}.

We organize this section as follows. After a brief description regarding the setting of the numerical simulations (Section~\ref{subsec: settings}), we compare the results of a physiological electromechanical simulation with a comprehensive set of experimental data available in literature (Section~\ref{subsec: baseline}). Then, in Section~\ref{subsec: cross_fibers} we investigate how different cross-fibers active contraction arrangements 
affect the electromechanical simulations, by setting different combinations of $n_{\text{f}}$, $n_{\text{s}}$ and $n_{\text{n}}$, i.e. of the prescribed proportion of active tension along the myofibers. Finally, in Section~\ref{subsec: configuration_fibers} we evaluate the impact of different myofiber architectures, obtained by three types of LDRBMs, on the biventricular pumping function.

\subsection{Setting of numerical simulations}
\label{subsec: settings}
All the simulations are performed on a realistic biventricular geometry processed from the Zygote 3D heart~\cite{zygote2014}, a CAD-model representing an average healthy human heart reconstructed from high-resolution computer tomography scan. To build the computational mesh associated with the biventricular Zygote model, we use the Vascular Modeling Toolkit software~\cite{antiga2008image} (\url{http://www.vmtk.org}) by exploiting the semi-automatic meshing tool recently proposed in~\cite{fedele2021polygonal}.

We employ two nested meshes where for the mechanical and activation problems we adopt a mesh size of $3\,\text{mm}$, while for the electrophysiology problem we employ a mesh size four time smaller~\cite{regazzoni2020numerical}.
As for the time steps, we use $\tau=50 \, \mu\text{s}$ for the
electrophysiology problem and $\Delta t=500 \, \mu\text{s}$ for the mechanical, activation and circulation problems~\cite{regazzoni2020numerical,piersanti2021modeling}.

The parameters used for the 3D-0D model are listed in Tables~\ref{tab:EMA_params} and~\ref{tab:C_params}. The settings related to LDRBMs, adopted for prescribing the fiber architectures, will be specified for each case reported in Sections.~\ref{subsec: baseline}$\,$--$\,$\ref{subsec: configuration_fibers}. 

To approach the limit cycle, we initialize all the numerical simulations, for the coupled 3D-0D model, following the procedure illustrated in Section~\ref{subsec: initial_conditions} (see also~\cite{regazzoni2021emulator}). Then, we perform three further heartbeats using the fully framework of the 3D-0D model presented in Sections~\ref{sec: mathematicalmodeling} and~\ref{sec: numericaldiscretization}. We neglected the first two, so that all the reported results refer to the last heartbeat.

In all the simulations we adopted the same pacing protocol in which five ventricular endocardial areas are activated with spherical impulses: in the anterior para-septal wall, in the left surface of inter-ventricular septum and in the bottom of postero-basal area, for LV; in the septum and in the free endocardial wall, for RV~\cite{piersanti2021modeling,durrer1970total}, see also Figure~\ref{fig: baseline}(c). This, combined with the fast endocardial conduction layer (see Section~\ref{subsec: electrophysiology}), surrogates the action of the Purkinje network~\cite{lee2019rule,vergara2014patient}.

The numerical methods presented in Section~\ref{sec: numericaldiscretization} have been implemented within \texttt{life\textsuperscript{x}} (\url{https://lifex.gitlab.io/lifex}), a new in-house high-performance \texttt{C++} FE library, for cardiac applications, based on the \texttt{deal.II} FE core \cite{dealII91}~(\url{https://www.dealii.org}). All the numerical simulations were executed using either the \texttt{iHeart} cluster (Lenovo SR950 192-Core Intel Xeon Platinum 8160, 2100 MHz and 1.7TB RAM) at MOX, Dipartimento di Matematica, Politecnico di Milano or the \texttt{GALILEO} supercomputer at Cineca (8~nodes endowed with 36 Intel Xeon E5-2697 v4 2.30 GHz).

\subsection{Baseline simulation}
\label{subsec: baseline}
\begin{figure}
	\centering
	\includegraphics[keepaspectratio, width=1\textwidth]{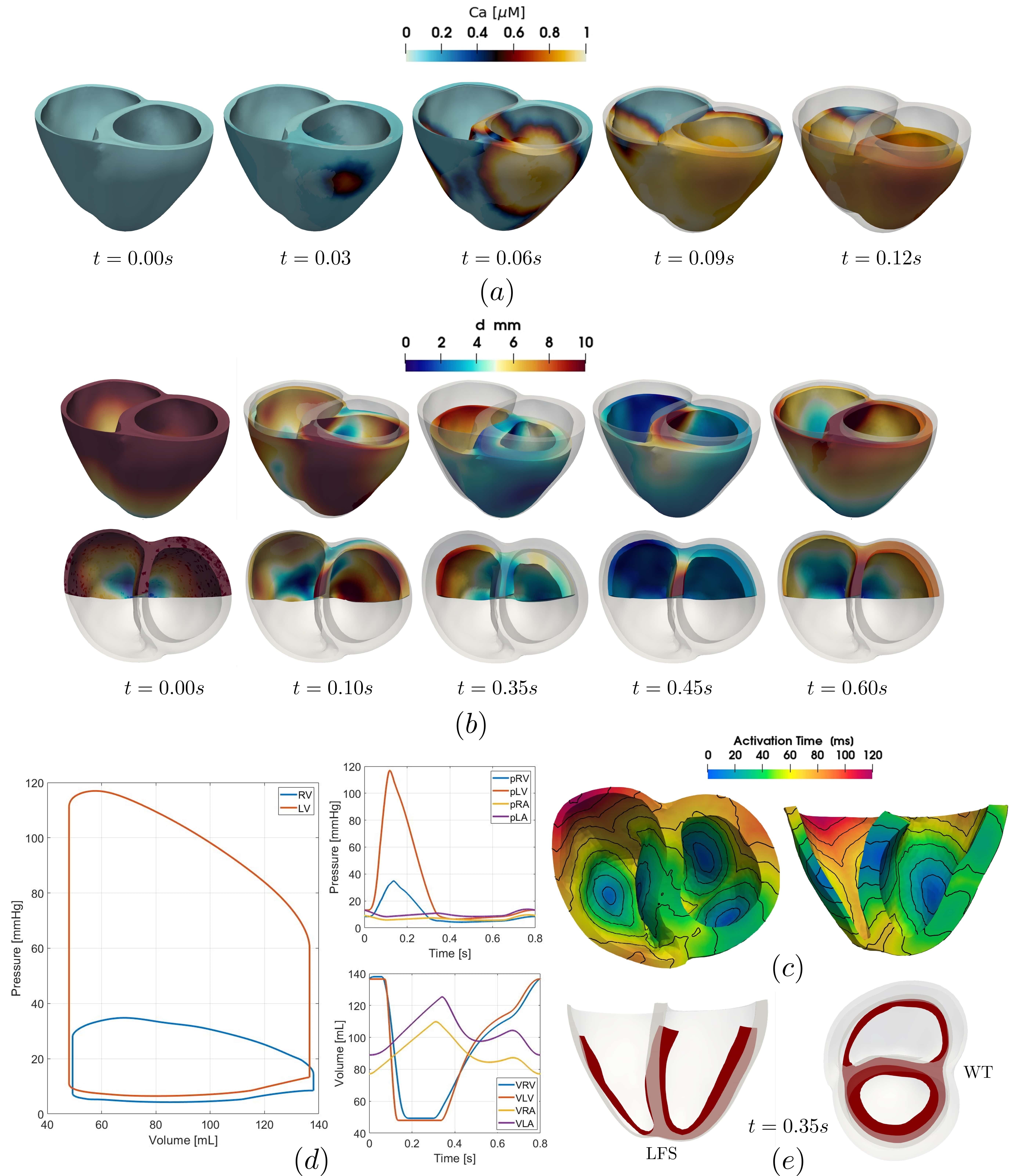}
	\caption{Baseline electromechanical simulation; (a) calcium transient at five time instants in the cardiac cycle; (b) mechanical displacement magnitude (with respect to the reference configuration) at five time instants of the heartbeat where $0.35\,s$ is the end of systole. (c) activation map; (d, left)  PV-loop LV (orange) and RV (blue); (d, right) pressures (top) and volumes (bottom) transient during the cardiac cycle for the four chambers; (e) mid ventricular slices at the end of systole, showing LFS on the left and WT on the right.}
	\label{fig: baseline}
\end{figure}
We present a human electromechanical simulation in physiological conditions related to the Zygote biventricular geometry. For the fibers generation we adopted D-RBM~\cite{piersanti2021modeling,doste2019rule}. The input angle values (see Section~\ref{subsec: fibers}) were chosen according to observations based on histological studies in the human heart~\cite{lombaert2012human,anderson2009three}:
\begin{equation}
	\label{eqn: angles}
	\begin{alignedat}{2}
		& \alpha_{epi,\text{LV}}=-60\degree, \quad \alpha_{endo,\text{LV}}=+60\degree, \quad \alpha_{epi,\text{RV}}=-25\degree, \quad \alpha_{endo,\text{RV}}=+90\degree;
		\\
		& \beta_{epi,\text{LV}}=+20\degree, \quad \beta_{endo,\text{LV}}=-20\degree, \quad \beta_{epi,\text{RV}}=+20\degree, \quad \beta_{endo,\text{RV}}=0\degree.
	\end{alignedat}
\end{equation}
Moreover, to surrogate the effect of dispersed myofibers, we set in~\eqref{eqn: piola} $n_{\text{f}}=0.7$, $n_{\text{s}}=0$ and $n_{\text{n}}=0.3$ for the proportion of active tension along the fiber, sheet and normal directions, respectively~\cite{guan2020effect,ahmad2018region}.

Figure~\ref{fig: baseline} illustrates the time evolution of calcium ions concentration (a), the mechanical deformation (b,~e), the activation times (c), the PV-loop curves for both ventricles and the time evolution of pressures and volumes of the four chambers (d). Specifically, in Figure~\ref{fig: baseline}(a) we display the time evolution of the TTP06 ionic model calcium transient showing the physiological wave propagation up to the complete depolarization of both ventricles ($t=0.12 \, \text{s}$). In Figure~\ref{fig: baseline}(b) we report different snapshots of the biventricular geometry warped by the displacement vector. As expected, at the beginning of the contraction the volumes of both ventricles remain nearly constant while the pressure increases ($t=0.0-0.10 \, \text{s}$); during the ejection phase, the ventricular contraction is clearly visible, with the basal plane that moves towards the bottom while the apex remains almost fixed. Moreover, a significant thickening of the myocardium wall takes place ($t=0.35 \, \text{s}$). Then, the ventricles start to relax. This leads to a slow recovery of the initial volumes ($t=0.45-0.60 \, \text{s}$). Finally, in Figure~\ref{fig: baseline}(c) we display the simulated activation map in which both the total activation time (120~ms) and the activation pattern are in accordance with the literature~\cite{piersanti2021modeling,durrer1970total}.

\begin{table}
	\centering
	\scalebox{1}{
		\begin{tabular}{lrrr}
			\toprule
			Mechanical biomarkers & Literature values & Simulation results & Description \\
			\midrule
			$\text{EDV}_{\text{LV}}$ (mL) & 142 $\pm$ 21 \cite{maceira2006normalized}  & 137 & Left end diastolic volume  \\
			$\text{EDV}_{\text{RV}}$ (mL) & 144 $\pm$ 23 \cite{maceira2006reference}  & 138 & Right end diastolic volume \\
			$\text{ESV}_{\text{LV}}$ (mL) & 47 $\pm$ 10 \cite{maceira2006normalized}  & 48 & Left end systolic volume  \\
			$\text{ESV}_{\text{RV}}$ (mL) & 50 $\pm$ 14 \cite{maceira2006reference} & 49 & Right end systolic volume  \\
			$\text{EF}_{\text{LV}}$ ($\%$) & 67 $\pm$ 5 \cite{maceira2006normalized}  & 66 & Left ventricular ejection fraction  \\
			$\text{EF}_{\text{RV}}$ ($\%$) & 67$\pm$ 8 \cite{tamborini2010reference}  & 65 & Right ventricular ejection fraction  \\
			$\text{P}_{\text{LV}}$ (mmHg) & 119 $\pm$ 13 \cite{sugimoto2017echocardiographic} & 117 & Left systolic pressure peak \\
			$\text{P}_{\text{RV}}$ (mmHg) & 35 $\pm$ 11 \cite{bishop1997clinical}  & 35 & Right systolic pressure peak  \\
			$\text{LFS}$ ($\%$) & 13-21 \cite{emilsson2006mitral}  & 21 & Longitudinal fractional shortening \\
			$\text{WT}$ ($\%$) & 18-100 \cite{sechtem1987regional}  & 41 & Fractional wall thickening \\
			\bottomrule
	\end{tabular}}
	\caption{Comparison between the simulation results, employing the 3D-0D EM model, and literature values of mechanical biomarkers in physiological conditions (references are reported in the Table, see also~\cite{levrero2020sensitivity,wang2021human}). }
	\label{tab:EM_comparison}
\end{table}

In Table~\ref{tab:EM_comparison} we compare some relevant mechanical biomarkers obtained by our numerical simulation with those provided by the data reported in the literature~\cite{maceira2006normalized,tamborini2010reference,maceira2006reference,sugimoto2017echocardiographic,bishop1997clinical,emilsson2006mitral,sechtem1987regional}. Notice that all the values in Table~\ref{tab:EM_comparison}, related to the ventricular volumes, are expressed with absolute values, in mL, estimated for an adult subject, as reported in the quoted references. However, we are aware that in the clinical practice the ventricular volumes are always indicated as "indexed ventricular volumes", by dividing the ventricular volume for the Body Surface Area of the related patient. The chosen mechanical biomarkers were: i) left and right end diastolic/systolic volumes ($\text{EDV}_{\text{LV}}$, $\text{EDV}_{\text{RV}}$, $\text{ESV}_{\text{LV}}$, $\text{ESV}_{\text{RV}}$), representing the maximal and minimal left and right ventricular volumes achieved during the heartbeat, computed as the maximal ($\text{EDV}_{\text{LV}}$, $\text{EDV}_{\text{RV}}$) and minimal ($\text{ESV}_{\text{LV}}$, $\text{ESV}_{\text{RV}}$) volumes in the PV-loop curves, see Figure~\ref{fig: baseline}(d); ii) left and right ventricular ejection fractions ($\text{EF}_{\text{LV}}$, $\text{EF}_{\text{RV}}$), which represent the amount of blood that is pumped by LV and RV during a cardiac cycle, computed as
$$\text{EF}_{\text{i}}=\frac{\text{EDV}_{\text{i}}-\text{ESV}_{\text{i}}}{\text{EDV}_{\text{i}}}100 \qquad \text{i}=\text{LV,RV};$$
iii) left and right systolic pressure peaks ($\text{P}_{\text{LV}}$, $\text{P}_{\text{RV}}$), the maximal pressures reached in LV and RV, computed as the maximal pressures in the PV-loop curves, see Figure~\ref{fig: baseline}(d); iv) the systolic longitudinal fractional shortening (LFS), standing for the fractional displacement between the endocardial apex and the base~\cite{levrero2020sensitivity}, evaluated as
$$\text{LFS}=\frac{\text{L}_{\text{0}}-\text{L}}{\text{L}_{\text{0}}}100,$$
where $\text{L}_{\text{0}}$ and $\text{L}$ are the apico-basal distance measured at the beginning ($t=0.0 \, \text{s}$) and at the end of systole ($t=0.35 \, \text{s}$), see Figure~\ref{fig: baseline}(e); v) the systolic wall thickening (WT), representing the fractional cardiac wall thickening~\cite{levrero2020sensitivity}, measured as
$$\text{WT}=\frac{\text{T}-\text{T}_{\text{0}}}{\text{T}}100,$$
where $\text{T}_{\text{0}}$ and $\text{T}$ are the cardiac wall thickening at the beginning ($t=0.0 \, \text{s}$) and at the end of systole ($t=0.35 \, \text{s}$), see Figure~\ref{fig: baseline}(e).

All the above mechanical biomarkers, obtained by our numerical simulation, fall within the physiological range (references in Table~\ref{tab:EM_comparison}).

\subsection{Cross-fibers active contraction}
\label{subsec: cross_fibers}
\begin{figure}
	\centering
	\includegraphics[keepaspectratio, width=1\textwidth]{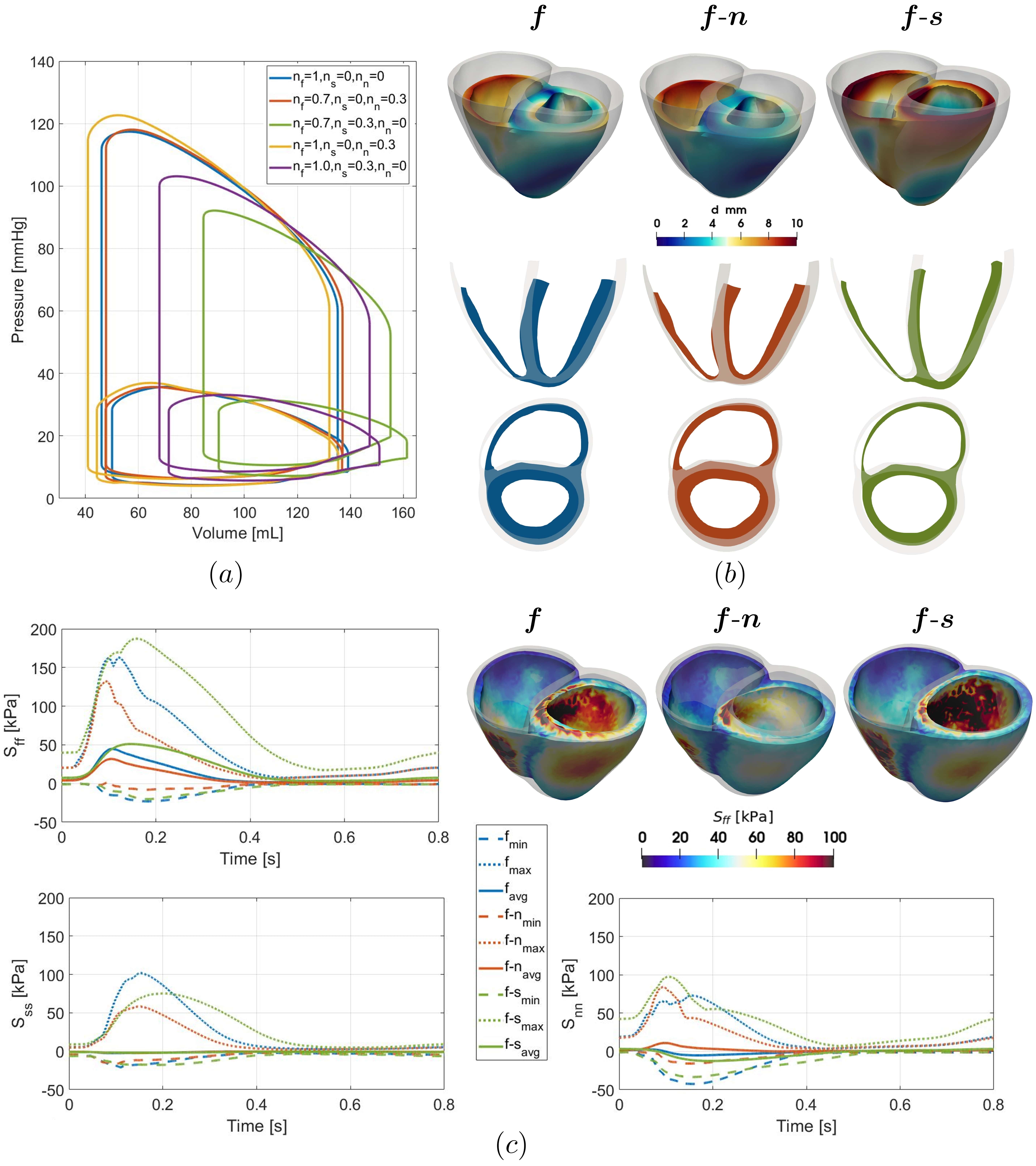}
	\caption{Cross-fibers active contraction simulations; (a) PV-loops from several cross-fibers active contraction arrangements built by setting in \eqref{eqn: piola} different combinations of $n_{\text{f}}$, $n_{\text{s}}$ and $n_{\text{n}}$; (b) mechanical displacements (top) and mid ventricular slices at the end of systole ($0.35\,s$), showing LFS (middle) and WT (bottom) for redistributed cross-fibers active contraction configurations: a pure fiber $\boldsymbol{f}$ (blue), a fiber-normal $\boldsymbol{f}$-$\boldsymbol{n}$ (orange) and a fiber-sheet $\boldsymbol{f}$-$\boldsymbol{s}$ (green) contractions; (c) circumferential stress $\stress_{\text{ff}}$ (top-right) at the peak pressure time instant ($0.1\,s$) and the time trace of the average, minimum and maximum axial stresses $\stress_{\text{ff}}$ (top-left), $\stress_{\text{ss}}$ (bottom-left) and $\stress_{\text{nn}}$ (bottom-right) for $\boldsymbol{f}$, $\boldsymbol{f}$-$\boldsymbol{n}$ and $\boldsymbol{f}$-$\boldsymbol{s}$ configurations.}
	\label{fig: cross_fibers}
\end{figure}
To surrogate the dispersion effect in the cardiac fibers, we analyse several cross-fibers active contraction
arrangements, by setting in \eqref{eqn: piola} different combinations of $n_{\text{f}}$, $n_{\text{s}}$ and $n_{\text{n}}$, i.e. the prescribed proportion of active tension along the myofibers. Five different sets were chosen: i) $n_{\text{f}}=0.7$, $n_{\text{s}}=0.3$, $n_{\text{n}}=0$; ii) $n_{\text{f}}=1$, $n_{\text{s}}=0.3$, $n_{\text{n}}=0$; iii) $n_{\text{f}}=1$, $n_{\text{s}}=0$, $n_{\text{n}}=0$; iv) $n_{\text{f}}=0.7$, $n_{\text{s}}=0$, $n_{\text{n}}=0.3$; v) $n_{\text{f}}=1$, $n_{\text{s}}=0$, $n_{\text{n}}=0.3$. Apart from the prescribed proportion of active tension, the settings are the same as the baseline simulation\footnote{Notice that case iv is the baseline simulation.} presented in Section~\ref{subsec: baseline}.

Figure~\ref{fig: cross_fibers}(a) shows the PV-loops from the five cases. An active tension along the sheet direction ($n_{\text{s}}>0$, cases i and ii) produces a PV-loop with a reduced area compared to case iii with no cross-fibers active contraction. Conversely, an active tension along the normal direction ($n_{\text{n}}>0$, cases iv and v) yields a PV-loop with an increased area. Table~\ref{tab:EF_cross-fibers} displays, for all the cases, the ejection fraction ($\text{EF}_{\text{i}}$) and the stroke volume ($\text{SV}_{\text{i}}=\text{EDV}_{\text{i}}-\text{ESV}_{\text{i}}$) of the
left ($\text{i}=\text{LV}$) and right ($\text{i}=\text{RV}$) ventricles. The maximal cardiac work is achieved for case v while the minimal for case i.
The above analysis shows that the active tension along the sheet direction ($n_{\text{s}}>0$) counteracts the myofiber contraction, while the one along the normal direction ($n_{\text{n}}>0$) enhances the cardiac work, in accordance to~\cite{guan2020effect,guan2021modelling}.
\begin{table}[t]
	\centering
	\scalebox{1}{
		\begin{tabular}{lrrrr}
			\toprule
			Cross-fiber configuration & $\text{EF}_{\text{LV}}$ & $\text{EF}_{\text{RV}}$ & $\text{SV}_{\text{LV}}$ & $\text{SV}_{\text{RV}}$ \\
			\midrule
			i) $n_{\text{f}}=0.7$, $n_{\text{s}}=0.3$, $n_{\text{n}}=0$ & 45$\,\%$ & 44$\,\%$ & 70.69 mL & 71.04 mL \\
			ii) $n_{\text{f}}=1$, $n_{\text{s}}=0.3$, $n_{\text{n}}=0$ & 54$\,\%$ & 53$\,\%$ & 79.40 mL & 79.50 mL  \\
			iii) $n_{\text{f}}=1$, $n_{\text{s}}=0$, $n_{\text{n}}=0$ & 65$\,\%$ & 64$\,\%$ & 89.14 mL & 89.08 mL \\
			iv) $n_{\text{f}}=0.7$, $n_{\text{s}}=0$, $n_{\text{n}}=0.3$ & 66$\,\%$ & 65$\,\%$ & 89.27 mL & 89.23 mL  \\
			v) $n_{\text{f}}=1$, $n_{\text{s}}=0$, $n_{\text{n}}=0.3$ & 69$\,\%$ & 67$\,\%$ & 91.14 mL & 91.09 mL \\
			\bottomrule
	\end{tabular}}
	\caption{Ejection fraction ($\text{EF}_{\text{i}}$) of the left ($\text{i}=\text{LV}$) and right ($\text{i}=\text{RV}$) ventricles for the different cross-fibers active contraction cases i$-$v. The stroke volume ($\text{SV}_{\text{i}}$) of the two ventricles is also shown.}
	\label{tab:EF_cross-fibers}
\end{table}

In order to better appreciate the differences among the cross-fibers active contraction arrangements, we further compared cases i and iv with case iii. In these particular cases, the proportion of active tension sums up to 1 ($n_\text{f}+n_\text{s}+n_\text{n}=1$), meaning that the myofibers contraction is redistributed along the three directions: case iii $(n_{\text{f}}=1$, $n_{\text{s}}=0$, $n_{\text{n}}=0)$ is a pure fiber contraction, in the following denoted by $\boldsymbol{f}$ configuration; case i ($n_{\text{f}}=0.7$, $n_{\text{s}}=0$, $n_{\text{n}}=0.3$) is a contraction in the fiber and normal directions,  hereafter indicated by $\boldsymbol{f}$-$\boldsymbol{n}$ configuration; case iv ($n_{\text{f}}=0.7$, $n_{\text{s}}=0.3$, $n_{\text{n}}=0$) is a contraction along the fiber and sheet directions, named $\boldsymbol{f}$-$\boldsymbol{s}$ configuration.

Figure~\ref{fig: cross_fibers}(b) illustrates the mechanical displacements at the end of systole ($0.35\,s$) for the three considered configurations ($\boldsymbol{f}$, $\boldsymbol{f}$-$\boldsymbol{n}$ and $\boldsymbol{f}$-$\boldsymbol{s}$). Both the apico-basal shortening and the wall thickening is dramatically reduced for $\boldsymbol{f}$-$\boldsymbol{s}$ configuration. Almost the same mechanical contraction is achieved for $\boldsymbol{f}$ and $\boldsymbol{f}$-$\boldsymbol{n}$ configurations with a slightly more pronounced longitudinal shortening and wall thickening for $\boldsymbol{f}$-$\boldsymbol{n}$ configuration. The LFS and WT are reported in Table~\ref{tab:LFT_WT}.
\begin{table}[t]
	\centering
	\scalebox{1}{
		\begin{tabular}{lrrr}
			\toprule
			Cross-fiber configuration & LFS & WT \\
			\midrule
			iii) $\boldsymbol{f}$ & 17$\,\%$ & 30$\,\%$ \\
			iv) $\boldsymbol{f}$-$\boldsymbol{n}$ & 21$\,\%$ & 41$\,\%$ \\
			i)  $ \: \: \boldsymbol{f}$-$\boldsymbol{s}$ & 7$\,\%$ & 8$\,\%$ \\
			\bottomrule
	\end{tabular}}
	\caption{LFS and WT for the three configurations of redistributed myofibers active contraction ($\boldsymbol{f}$, $\boldsymbol{f}$-$\boldsymbol{n}$ and $\boldsymbol{f}$-$\boldsymbol{s}$)}
	\label{tab:LFT_WT}
\end{table}

We also evaluate the components of the mechanical stress by means of the following indicators~\cite{regazzoni2020numerical}:
\begin{equation*}
	\stress_{\text{ff}} = (\tenspiola \mathbf{f}_{\text{0}}) \cdot \dfrac{\mecF \mathbf{f}_{\text{0}}}{|\mecF \mathbf{f}_{\text{0}}|}, \quad 	\stress_{\text{ss}} = (\tenspiola \mathbf{s}_{\text{0}}) \cdot \dfrac{\mecF \mathbf{s}_{\text{0}}}{|\mecF \mathbf{s}_{\text{0}}|}, \quad 	\stress_{\text{nn}} = (\tenspiola \mathbf{n}_{\text{0}}) \cdot \dfrac{\mecF \mathbf{n}_{\text{0}}}{|\mecF \mathbf{n}_{\text{0}}|},
\end{equation*}
where $\mathbf{f}_{\text{0}}$, $\mathbf{s}_{\text{0}}$ and $\mathbf{n}_{\text{0}}$ are the myofiber directions, $\tenspiola$ is the first Piola-Kirchhoff stress tensor and $\mecF$ is the deformation gradient tensor. The metric $\stress_{\text{aa}}$ (with $\text{a = f, s, n}$) measures the axial stresses along the circumferential ($\text{a = f}$), radial ($\text{a = s}$) and longitudinal ($\text{a=n}$) directions.

Figure~\ref{fig: cross_fibers}(c) displays, for the three configurations $\boldsymbol{f}$, $\boldsymbol{f}$-$\boldsymbol{n}$ and $\boldsymbol{f}$-$\boldsymbol{s}$,  the circumferential stress ($\stress_{\text{ff}}$) at the peak pressure time instant ($0.1\,s$) and the time trace of the average, minimum and maximum axial stresses $\stress_{\text{ff}}$, $\stress_{\text{ss}}$ and $\stress_{\text{nn}}$. The circumferential stress at the peak pressure instant is much higher, especially on LV side, for $\boldsymbol{f}$-$\boldsymbol{s}$ configuration with respect to the other two. Conversely, $\boldsymbol{f}$-$\boldsymbol{n}$ configuration produces the lowest circumferential stress. Almost the same considerations hold for the time trace of the three axial stresses during the complete cardiac cycle, see Figure~\ref{fig: cross_fibers}(c).

The previous results reveal that the configuration $\boldsymbol{f}$-$\boldsymbol{n}$ allows to obtain a more efficient cardiac contraction with a much lower axial stress with respect to $\boldsymbol{f}$ configuration. On the contrary, $\boldsymbol{f}$-$\boldsymbol{s}$ configuration yields an unphysiological cardiac contraction with EF, LFS and WT
below the physiological range reported in literature (see Tables~\ref{tab:EM_comparison}-\ref{tab:LFT_WT}).

\subsection{Impact of myofiber architecture on the electromechanical function}
\label{subsec: configuration_fibers}
\begin{figure}
	\centering
	\includegraphics[keepaspectratio, width=0.95\textwidth]{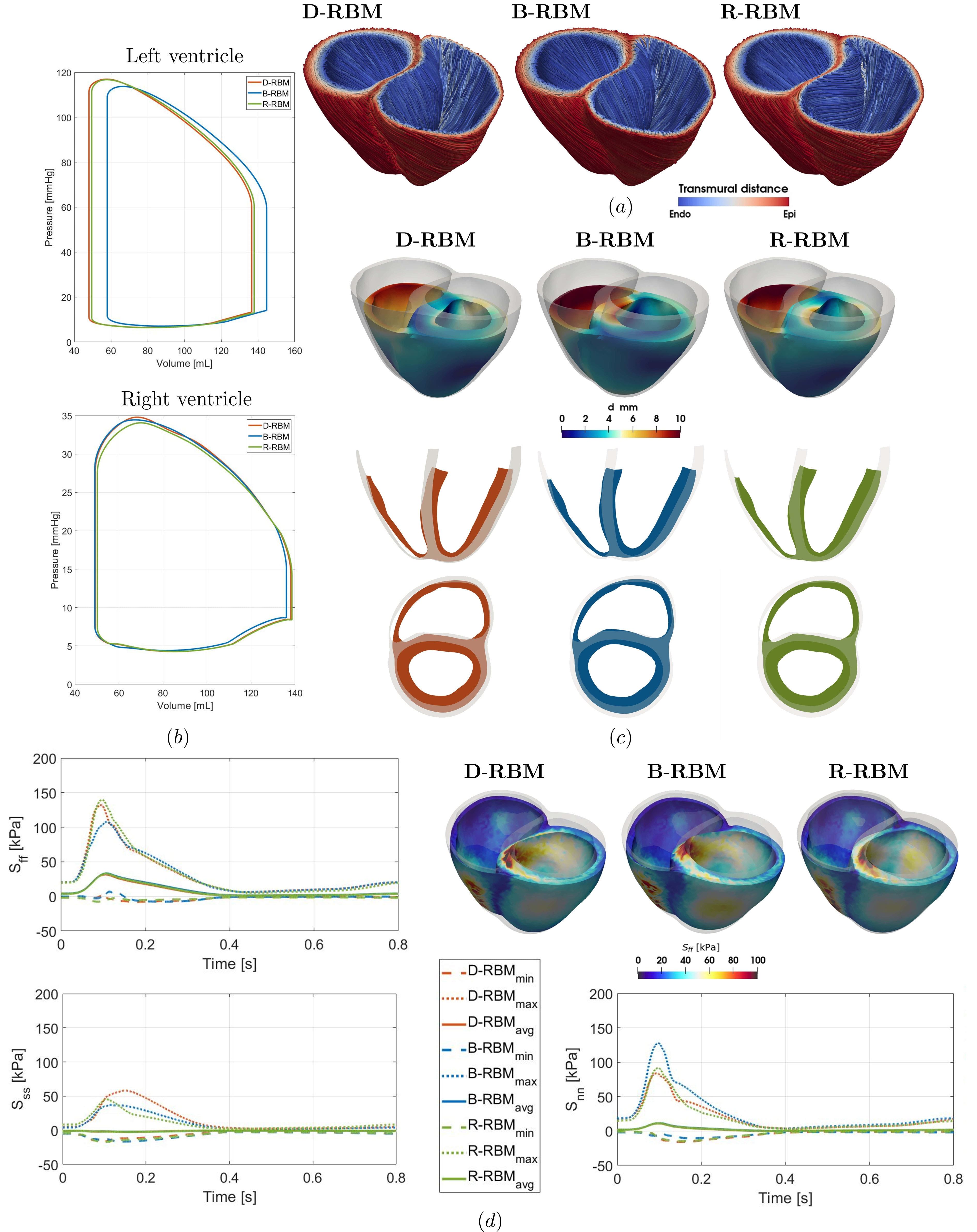}
	\caption{Results of the EM model employing different LDRBMs (R-RBM, B-RBM and D-RBM) to generate the fiber architecture; (a) fiber orientations obtained for the three LDRBMs in the Zygote biventricular model; (b)  PV-loop curves, for LV (top) and RV (bottom), obtained with the three LDRBMs: D-RBM (orange), B-RBM (blue) and R-RBM (green); (c) mechanical displacements (top) and mid ventricular slices at the end of systole ($0.35\,s$), showing LFS (middle) and WT (bottom) obtained by D-RBM (orange), B-RBM (blue) and R-RBM (green); (d) circumferential stress $\stress_{\text{ff}}$ (top-right) at the peak pressure instant ($0.1\,s$) and the time trace of the average, minimum and maximum axial stresses $\stress_{\text{ff}}$ (top-left), $\stress_{\text{ss}}$ (bottom-left) and $\stress_{\text{nn}}$ (bottom-right) for the three LDRBMs.}
	\label{fig: RBM_fibers}
\end{figure}
We investigate the effect of different myofibers architecture on the biventricular EM model, by considering three types of LDRBMs: D-RBM, B-RBM and R-RBM (see Section~\ref{subsec: fibers}). Apart from the employed LDRBM, used to prescribe the myofibers architecture, all the other settings, including the fiber input angles~\eqref{eqn: angles}, are the same as the baseline simulation\footnote{Notice that the case with D-RBM is the baseline simulation.} presented in Section~\ref{subsec: baseline}.

Fiber orientations obtained for the three LDRBMs (D-RBM, B-RBM and R-RBM) in the Zygote biventricular model are shown in Figure~\ref{fig: RBM_fibers}(a). For a detailed comparison among the three LDRBMs we refer the reader to~\cite{piersanti2021modeling}, where pure electrophysiological simulations were considered. Here, we
are instead interested in the effect on mechanical quantities obtained by means of EM model. We recall that B-RBM produces a smooth change in the fiber field in the transition across the two ventricles, while R-RBM and D-RBM a strong discontinuity~\cite{piersanti2021modeling}. Moreover, R-RBM and D-RBM feature a linear transition passing from the endocardium to the epicardium, while B-RBM employs a bidirectional spherical interpolation~\textit{bislerp} (see~\cite{piersanti2021modeling,bayer2012novel,doste2019rule,quarteroni2017integrated}).

In Figure~\ref{fig: RBM_fibers}(b) the PV-loop curves (for both ventricles) are displayed, while in Table~\ref{tab:RBM_comparison} some relevant mechanical biomarkers are compared among the simulation results. The left ventricular PV-loop area of B-RBM is shifted towards larger volumes with respect to the ones of D-RBM and R-RBM that show almost a compatible PV-loop for LV, see Figure~\ref{fig: RBM_fibers}(b,~top). Moreover, the left systolic pressure peak decreases for B-RBM with respect to D-RBM and R-RBM, see Figures~\ref{fig: RBM_fibers}(b, top) and Table~\ref{tab:RBM_comparison}. As a consequence, the left ventricular ejection fraction obtained with B-RBM ($60\%$) is smaller than those obtained with D-RBM and R-RBM ($66\%$ and $65\%$, respectively), see Table~\ref{tab:RBM_comparison}. On the contrary, small differences are observed for the right ventricular PV-loops with only a slightly larger ejection fraction for D-RBM, see Figure~\ref{fig: RBM_fibers}(b, bottom) and Table~\ref{tab:RBM_comparison}.
\begin{table}
	\centering
	\scalebox{1}{
		\begin{tabular}{lrrr}
			\toprule
			Mechanical biomarkers & D-RBM & B-RBM & R-RBM \\
			\midrule
			$\text{EDV}_{\text{LV}}$ (mL) & 137 & 145 & 138 \\
			$\text{EDV}_{\text{RV}}$ (mL) & 138 & 136 & 139 \\
			$\text{ESV}_{\text{LV}}$ (mL) & 48 & 58 & 50 \\
			$\text{ESV}_{\text{RV}}$ (mL) & 49 & 49 & 50 \\
			$\text{EF}_{\text{LV}}$ ($\%$) & 66 & 60 & 64 \\
			$\text{EF}_{\text{RV}}$ ($\%$) & 65 & 64 & 64 \\
			$\text{P}_{\text{LV}}$ (mmHg) & 117 & 114 & 117  \\
			$\text{P}_{\text{RV}}$ (mmHg) & 35 & 34 & 33  \\
			$\text{LFS}$ ($\%$) & 21 & 25 & 20 \\
			$\text{WT}$ ($\%$) & 41 & 36 & 38 \\
			\bottomrule
	\end{tabular}}
	\caption{Comparison of relevant mechanical biomarkers among the electromechanical simulations by employing different LDRBMs (D-RBM,B-RBM and R-RBM) to prescribe the myofiber architecture.}
	\label{tab:RBM_comparison}
\end{table}

Figure~\ref{fig: RBM_fibers}(d) shows the circumferential stress ($\stress_{\text{ff}}$) at the peak pressure instant ($0.1\,s$) and the time trace of the average, minimum and maximum axial stresses $\stress_{\text{ff}}$, $\stress_{\text{ss}}$ and $\stress_{\text{nn}}$. The patterns of $\stress_{\text{ff}}$ are very similar for the three methods, see Figure~\ref{fig: RBM_fibers}(d, top-right). Instead, the time traces of the axial stresses present several discrepancies. Specifically, $\stress_{\text{ff}}$ reveals lower values obtained by B-RBM with respect to D-RBM and R-RBM, see Figure~\ref{fig: RBM_fibers}(d, top-left). This is associated to a lower cardiac work produced by B-RBM ($\text{EF}_{\text{LV}}=60\%$) compared to D-RBM and R-RBM ($\text{EF}_{\text{LV}}=66\%$, $64\%$, respectively). On the contrary, the longitudinal stress $\stress_{\text{nn}}$ presents an opposite trend, see Figure~\ref{fig: RBM_fibers}(d, bottom-right). This is ascribed to a larger apico-basal shortening for B-RBM ($\text{LFS}=25\%$) with respect to D-RBM and R-RBM ($\text{LFS}=21\%$, $20\%$, respectively). Meanwhile, larger values of the radial stress $\stress_{\text{ss}}$ are observed for D-RBM with respect to B-RBM and R-RBM, see Figure~\ref{fig: RBM_fibers}(d,bottom-right), associated to a larger wall thickening of D-RBM ($\text{WT}=41\%$) against the ones of R-RBM and B-RBM ($\text{WT}=38\%$, $36\%$, respectively).

The previous results highlight that there is a strong interaction on the cardiac pump function between the LV and RV~\cite{palit2015computational}. A different fibers architecture in the transmural wall (from epicardium to endocardium) and a different septal fibers interconnection between the two ventricles affect the ventricular cardiac pump work, in particular the LV one. Indeed, a biventricular myofibers architecture has much more information (e.g. in the inter-ventricular septum) compared to a stand-alone LV model.

%% file: parts_conclusions.tex
\section{Conclusions}
\label{sec: conclusions}
In this work, we presented a 3D biventricular EM model coupled with a 0D closed-loop model of the whole cardiovascular system. We provided a rigorous mathematical and numerical formulation 
of the 3D-0D model by fully detailing our approach to couple the 3D and the 0D models.
We carried out several numerical simulations aimed at reproducing physiological quantities
like the PV-loops. Our results quantitatively match the experimental data of relevant mechanical biomarkers available in literature~\cite{maceira2006normalized,tamborini2010reference,maceira2006reference,sugimoto2017echocardiographic,bishop1997clinical,emilsson2006mitral,sechtem1987regional}, such as the end systolic and diastolic volumes, the ejection fractions, the systolic pressure peaks, the longitudinal fractional shortening and the fractional wall thickening. 

We studied different configurations in cross-fibers active contraction proving that an active tension along the sheet-normal direction enhances the cardiac work, whereas along the sheet direction it has the opposite effect.
Moreover, an active contraction in the sheet-normal direction allows to obtain a more efficient cardiac pumping function with a much lower axial stress with respect to a pure fiber configuration. Conversely, a sheet active contraction yields unphysiological ejection fraction, longitudinal shortening and wall thickening. These results put in evidence that the proportion of active tension along the sheet direction should be avoided in the framework of an orthotropic active stress.

Finally, we evaluate the impact of different myofibers architecture on the biventricular EM. Our results showed the importance of considering a biventricular model with respect to a stand-alone LV model. A different fibers architecture in the transmural wall and in the inter-ventricular septum influence the ventricular cardiac pump work, in particular the LV one. This highlights the strong interaction on the cardiac pump function between the LV and RV, highlighting the importance of considering the two chambers together during the ventricular electromechanical simulation. The continuous interrelationships between right and left ventricular functions are well known not only in physiological conditions, but particularly in pathological situations, for which any pressure and/or volume overload of a ventricle is instantaneously reflected in impairment of the function of the contralateral ventricle.

%% file: acknowledgements.tex
\section*{Acknowledgements}
This project has received funding from the European Research Council (ERC) under the European Union's Horizon 2020 research and innovation programme (grant agreement No 740132, iHEART - An Integrated Heart Model for the simulation of the cardiac function, P.I. Prof. A. Quarteroni). We acknowledge the CINECA award under the class C ISCRA project HP10C3Z520 for the availability of high performance computing resources.  
\begin{center}
	\raisebox{-.5\height}{\includegraphics[width=.15\textwidth]{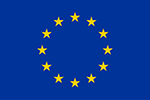}}
	\hspace{2cm}
	\raisebox{-.5\height}{\includegraphics[width=.15\textwidth]{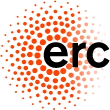}}
\end{center}

%% file: parts_app_params.tex
\section{Model parameters}
\label{app:params}
We provide the list of parameters adopted for the simulations in Sec.~\ref{sec: numericalresults}.
In particular, Table~\ref{tab:EMA_params} contains the parameters of the 3D EM model (referred to $\PhyEP$, $\PhyAct$, $\PhyMec$) and Table~\ref{tab:C_params} those of the 0D closed-loop hemodynamical model ($\PhyCirc$).
Moreover, for the TTP06 ionic model, we use the parameters (for epicardium cells) reported in~\cite{ten2006alternans}, while for the RDQ18 model, we employ those in~\cite{regazzoni2020machine}. 
\begin{table}[ht!]
	\centering
	\scalebox{1}{
	\begin{tabular}{lrrr}
		\toprule
		Variable & Value & Unit & Description \\
		\midrule
		\multicolumn{4}{l}{\textbf{Electrophysiology}} \\
		$T_{\text{hb}}$ & \num{0.8} & $\text{s}$ & Heartbeat duration \\
		$\EPchim$ & \num{1} & $\mu \text{F}/\text{cm}^2$ & Surface-to-volume ratio \\
		$\EPCm$ & \num{1400} & $\text{cm}^{-1}$& Transmembrane
		capacitance \\
		$\epsilon$ & \num{0.01} & $-$ & Threshold of the fast conduction layer \\
		$(\sigma_{\ell,\text{fast}}, \sigma_{\text{t,fast}}, \sigma_{\text{n,fast}})$ & (\num{4.28}, \num{1.96}, \num{0.64}) & $\text{mS/cm}$ & Fast layer conductivities \\
		$(\sigma_{\ell,\text{myo}}, \sigma_{\text{t,myo}}, \sigma_{\text{n,myo}})$ & (\num{1.07}, \num{0.49}, \num{0.16}) & $\text{mS/cm}$ & Myocardial conductivities \\
	    $\EPIapp$ & \num{50e3} & $\mu \text{A}/\text{cm}^3$ & Applied current value \\
	    $\tIapp$ & \num{3.0} & ms & Applied current duration \\
	    $\tIappzLV$ & (\num{0.0},\num{0.0},\num{0.0}) & ms & Applied current LV initial times \\
	    $\tIappzRV$ & (\num{5.0},\num{5.0}) & ms & Applied current RV initial times \\
	    $\DIapp$ & \num{2.5e-3} & \si{\meter} & Applied current radius \\
		\toprule
		\multicolumn{4}{l}{\textbf{Mechanics}} \\
		$\rho_\text{s}$  & $10^3$       & \si{\kilogram \per \cubic\meter} & Tissue density \\
	    $\BCmecKepiT$    & \num{2e4}    & \si{\pascal\per\meter}        & Normal stiffness of epicardium \\
	    $\BCmecKepiN$    & \num{2e5}    & \si{\pascal\per\meter}        & Tangential stiffness of epicardial tissue  \\
		$\BCmecCepiN$    & \num{2e4}    & \si{\pascal\second\per\meter} & Normal viscosity of epicardial tissue  \\
		$\BCmecCepiT$    & \num{2e3}    & \si{\pascal\second\per\meter} & Tangential viscosity of epicardial tissue \\
		$a$              & \num{0.88e3} & \si{\pascal}                  & Material stiffness \\
		$k$              & \num{50e3}   & \si{\pascal}                  & Bulk modulus \\
		$b_{\text{ff}}$  & 8            & $-$                           & Fiber strain scaling \\
		$b_{\text{ss}}$  & 6            & $-$                           & Radial strain scaling \\
		$b_{\text{nn}}$  & 3            & $-$                           & Cross-fiber in-plain strain scaling \\
		$b_{\text{fs}}$  & 12           & $-$                           & Shear strain in fiber-sheet plane scaling \\
	    $b_{\text{fn}}$  & 3            & $-$                           & Shear strain in fiber-normal plane scaling \\
	    $b_{\text{sn}}$  & 3            & $-$                           & Shear strain in sheet-normal plane scaling \\
	    \toprule
	    \multicolumn{4}{l}{\textbf{Reference Configuration}} \\
	    $\pressLVRef$          & 600            & \si{\pascal}                & Residual left ventricular pressure \\
	    $\pressRVRef$          & 400            & \si{\pascal}                & Residual right ventricular pressure \\
	    $\tensRef$       & \num{350e3}  & \si{\pascal}                     & Residual active tension \\
	    $\Cr$            & 1            & $-$                           & Residual contractility ratio \\
	    \toprule
		\multicolumn{4}{l}{\textbf{Activation}} \\
	    $\SL_0$          & 2            & \si{\micro\meter}                & Reference sarcomere length \\
	    $\TensMax$       & \num{840e3}  & \si{\pascal}                     & Maximum tension \\
	    $\Cr$            & 0.60            & $-$                           & Contractility ratio \\
		\bottomrule
	\end{tabular}}
	\caption{Input parameters of the 3D EM model.}
	\label{tab:EMA_params}
\end{table}
\begin{table}[ht!]
	\centering
	\scalebox{1}{
	\begin{tabular}{lrrr}
		\toprule
		Variable & Value & Unit & Description \\
		\midrule
		\multicolumn{4}{l}{\textbf{Circulation}} \\
		$\RarSYS$    & 0.416         & \si{\mmHg \second \per \milli\liter}         & Resistance of systemic arterial system \\
	    $\RvnSYS$    & 0.260        & \si{\mmHg \second \per \milli\liter}         & Resistance of systemic venous system \\
	    $\RarPUL$    & 0.048      & \si{\mmHg \second \per \milli\liter}         & Resistance of pulmonary arterial system \\
	    $\RvnPUL$    & 0.036      & \si{\mmHg \second \per \milli\liter}         & Resistance of pulmonary venous system \\
		$\CarSYS$    & 1.62         & \si{\milli\liter \per \mmHg}                 & Capacitance of systemic arterial system \\
		$\CvnSYS$    & 60.00        & \si{\milli\liter \per \mmHg}                 & Capacitance of systemic venous system \\
		$\CarPUL$    & 5.00        & \si{\milli\liter \per \mmHg}                 & Capacitance pulmonary arterial system \\
		$\CvnPUL$    & 16.00        & \si{\milli\liter \per \mmHg}                 & Capacitance of pulmonary venous system \\
		$\LarSYS$    & \num{5e-3}  & \si{\mmHg \second\squared \per \milli\liter} & Impedance of systemic arterial system \\
		$\LvnSYS$    & \num{5e-4}  & \si{\mmHg \second\squared \per \milli\liter} & Impedance of systemic venous system \\
		$\LarPUL$    & \num{5e-4}  & \si{\mmHg \second\squared \per \milli\liter} & Impedance pulmonary arterial system \\
		$\LvnPUL$    & \num{5e-4}  & \si{\mmHg \second\squared \per \milli\liter} & Impedance of pulmonary venous system \\
		$\EpLA$      & 0.09        & \si{\mmHg \per \milli\liter}                 & Left atrium elastance amplitude \\
		$\EpRA$      & 0.06        & \si{\mmHg \per \milli\liter}                 & Right atrium elastance amplitude  \\
		$\EaLA$      & 0.07        & \si{\mmHg \per \milli\liter}                 & Left atrium elastance baseline \\
		$\EaRA$      & 0.07        & \si{\mmHg \per \milli\liter}                 & Right atrium elastance baseline \\
        $\TacLA$     & 0.17        & $-$                                  & Duration of left atrium contraction (w.r.t. $T_{\text{hb}}$) \\
        $\TacRA$     & 0.17         & $-$                                 & Duration of right atrium contraction (w.r.t. $T_{\text{hb}}$) \\
        $\tacLA$     & 0.80         & $-$                                 & Initial time of left atrium contraction (w.r.t. $T_{\text{hb}}$) \\
        $\tacRA$     & 0.80         & $-$                                 & Initial time of right atrium  contraction (w.r.t. $T_{\text{hb}}$) \\
		$\TarLA$     & 0.17        & $-$                                 & Duration of left atrium relaxation (w.r.t. $T_{\text{hb}}$) \\
		$\TarRA$     & 0.17         & $-$                                 & Duration of right atrium relaxation (w.r.t. $T_{\text{hb}}$) \\
		$\VnLA$      & 4.0         & \si{\milli\liter}                            & Left atrium resting volume \\
		$\VnRA$      & 4.0         & \si{\milli\liter}                            & Right atrium resting volume \\
		$\Rmin$      & \num{75e-4}      & \si{\mmHg \second \per \milli\liter}         & Valves minimal resistance \\
		$\Rmax$      & \num{75e3}     & \si{\mmHg \second \per \milli\liter}         & Valves maximum resistance \\
		\bottomrule
	\end{tabular}}
	\caption{Input parameters of the 0D closed-loop hemodynamical model.}
	\label{tab:C_params}
\end{table}

%% file: parts_app_BC.tex
\section{Energy-consistent boundary condition in biventricular geometries}
\label{app:BC}

The energy-consistent boundary condition~\eqref{eqn: Mec5} accounts for the effect of the neglected part of the domain located above the biventricular base $\GammaBase$ (which is an artificial boundary), consistently with the principles of momentum and energy conservation.
It represents a generalization of the boundary condition proposed in \cite{regazzoni2020machine} for biventricular geometries.
In what follows, we denote by $\OmegaFluidtL$ (respectively $\OmegaFluidtR$) the volume occupied at time $t$, within LV (respectively, RV), by the fluid located below the base.
Moreover, we employ the tilde symbol ($\sim$) to refer to volumes and surfaces located above the ventricular base.
Specifically, we denote by $\OmegaFluidtTildeL$ and $\OmegaFluidtTildeR$ the fluid volumes in LV and RV, located above the base.
Similarly, we denote by $\GammaEpitTilde$,
$\GammaEndotTildeL$ and
$\GammaEndotTildeR$ the epicardial, and endocardial (left and right) surfaces located above the ventricular base.
Finally, we denote by $\GammaBasetTilde$ the ventricular base surface itself, but endowed with outer normal vector directed towards the apex, differently than for $\GammaBaset$.

Following the derivation of \cite{regazzoni2020machine} and by defining the Cauchy stress tensor as $\tenscaychy = J^{-1} \tenspiola \mecF^T$, with a quasi-static approximation the balance of momentum entails
\begin{equation}\label{eqn:BC_summation}
	\begin{split}
	\mathbf{0} &= \int_{\widetilde{\Omega}_t} \nabla \cdot \tenscaychy \, d\mathbf{x} =
	\int_{\partial\widetilde{\Omega}_t} \tenscaychy\mecNact \, d\Gamma_t
	=\int_{\GammaEpitTilde} \tenscaychy\mecNact \, d\Gamma_t +
	\int_{\GammaEndotTildeL} \tenscaychy\mecNact \, d\Gamma_t +
	\int_{\GammaEndotTildeR} \tenscaychy\mecNact \, d\Gamma_t +
	\int_{\GammaBasetTilde} \tenscaychy\mecNact \, d\Gamma_t.
	\end{split}
\end{equation}
The normal stress on the endocardium is given by $\tenscaychy\mecNact = -\PLV \mecNact$ (on $\GammaEndotTildeL$) and $\tenscaychy\mecNact = -\PRV \mecNact$ (on $\GammaEndotTildeR$), while we assume negligible the load on the epicardium (i.e. $\tenscaychy\mecNact = \mathbf{0}$ on $\GammaEpitTilde$).
Thanks to the divergence (Gauss) theorem, it is possible to write the endocardial terms of the summation of Eq.~\eqref{eqn:BC_summation} as integrals over $\GammaEndotL$ and
$\GammaEndotR$.
Indeed, we have the identity:
\begin{equation*}
\mathbf{0}
= \int_{\OmegaFluidtL \cup \OmegaFluidtTildeL} \nabla \PLV \, d\mathbf{x}
= \int_{\GammaEndotL} \PLV\mecNact \, d\Gamma_t
+ \int_{\GammaEndotTildeL} \PLV\mecNact \, d\Gamma_t,
\end{equation*}
and similarly for the RV we have $\textstyle\int_{\GammaEndotTildeL} \PLV\mecNact \, d\Gamma_t
= - \textstyle\int_{\GammaEndotL} \PRV\mecNact \, d\Gamma_t$.
Hence, we end up with the following identity
\begin{equation*}
	\begin{split}
	\int_{\GammaBaset} \tenscaychy\mecNact \, d\Gamma_t
	&= - \int_{\GammaBasetTilde} \tenscaychy\mecNact \, d\Gamma_t
	=
	-\int_{\GammaEndotTildeL} \PLV \mecNact \, d\Gamma_t
	-\int_{\GammaEndotTildeR} \PRV \mecNact \, d\Gamma_t
	\\&=
	\int_{\GammaEndotL} \PLV \mecNact \, d\Gamma_t +
	\int_{\GammaEndotR} \PRV \mecNact \, d\Gamma_t,
	\end{split}
\end{equation*}
which entails, by considering the pull-back, to the reference configuration
\begin{equation}\label{eqn:BC_integral}
	\begin{split}
	\int_{\GammaBaset} \tenscaychy\mecNact \, d\Gamma_t
	&=
	\int_{\GammaEndoL} \PLV \,J \mecF^{-T}\mecNact \, d\Gamma_0 +
	\int_{\GammaEndoR} \PRV \,J \mecF^{-T}\mecNact \, d\Gamma_0.
	\end{split}
\end{equation}
Equation~\eqref{eqn:BC_integral} provides the overall stress acting on the ventricular base.
However, we need some additional assumptions to define the point-wise distribution of stress, among the infinitely many satisfying Eq.~\eqref{eqn:BC_integral}.
In the original derivation of the energy-consistent boundary condition \cite{regazzoni2020machine}, at this stage, a uniform stress distribution assumption is made.
However, while this assumption is reasonable in a single-ventricle geometry, it is unrealistic when the ventricular base surrounds both ventricles.
Indeed, the pressures acting in LV are typically much larger than those in RV.
For this reason, we propose to distribute stress over the gamma surface not uniformly, but rather according to a weight function $\BCweight \colon \GammaBase \to [0,1]$, that indicates the fraction of stress attributable to the pressure acting on LV, relative RV, at each base point.
Hence, we assume that, on $\GammaBase$, we have:
\begin{equation}
	\begin{split}
		\tenscaychy\mecNact
		&=
		\BCweight
		\frac{\int_{\GammaEndoL} \PLV \,J \mecF^{-T}\mecNact \, d\Gamma_0}
		{\int_{\GammaBaset} \BCweight \, d\Gamma}
		+
		(1-\BCweight)
		\frac{\int_{\GammaEndoR} \PRV \,J \mecF^{-T}\mecNact \, d\Gamma_0}
		{\int_{\GammaBaset} (1-\BCweight) \, d\Gamma},
	\end{split}
\end{equation}
which reads, in the reference configuration:
\begin{equation} \label{eqn:BC_generic}
	\begin{split}
		\tenspiola\mecNact
		=
		{|J \mecF^{-T} \mecNact|}
		&\left[
			\BCweight
		\frac{\int_{\GammaEndoL} \PLV \,J \mecF^{-T}\mecNact \, d\Gamma_0}
		{\int_{\GammaBase} |J \mecF^{-T} \mecNact| \, \BCweight \, d\Gamma_0}
		+
		(1-\BCweight)
		\frac{\int_{\GammaEndoR} \PRV \,J \mecF^{-T}\mecNact \, d\Gamma_0}
		{\int_{\GammaBase} |J \mecF^{-T} \mecNact| \, (1-\BCweight)  \, d\Gamma_0}
		\right].
	\end{split}
\end{equation}
In what follows we consider three different choices for the weight function $\BCweight$, corresponding to as many boundary condition formulations.

\begin{itemize}
	\item \textbf{Uniform stress distribution}.
	By setting $\BCweight \equiv \tfrac{1}{2}$, we recover the case of stress uniformly distributed on the whole $\GammaBase$ boundary:
	\begin{equation}
	\begin{split}
	\tenspiola\mecNact &=
	\frac{|J \mecF^{-T} \mecNact|}
	{\int_{\GammaBase} |J \mecF^{-T} \mecNact| d\Gamma_0 }
	\left[ \int_{\GammaEndoL} \PLV \,J \mecF^{-T}\mecNact \, d\Gamma_0 + \right. \\
	& \qquad \qquad \qquad \qquad \quad \left. + \int_{\GammaEndoR} \PRV \,J \mecF^{-T}\mecNact \, d\Gamma_0 \right]
	\end{split}
	\end{equation}

	\item \textbf{Uniform stress distribution over each base}.
	Let us suppose to split the base into two subsets $\GammaBaseL$ and $\GammaBaseR$, respectively denoting the portion of ventricular base surrounding LV and RV.
	Then, we define $\BCweight$ as the indicator function of the set $\GammaBaseL$ (that is $\BCweight = 1$ on $\GammaBaseL$, while $\BCweight = 0$ on $\GammaBaseR$).
	In this case, we get:
	\begin{equation}
	\left\{
	\begin{split}
	\tenspiola\mecNact &= \frac{|J \mecF^{-T} \mecNact|}
	{\int_{\GammaBaseL} |J \mecF^{-T} \mecNact| d\Gamma_0} \int_{\GammaEndoL} \PLV \,J \mecF^{-T}\mecNact \, d\Gamma_0
	\qquad \text{on $\GammaBaseL$} \\
	\tenspiola\mecNact &= \frac{|J \mecF^{-T} \mecNact|}
	{\int_{\GammaBaseR} |J \mecF^{-T} \mecNact| d\Gamma_0} \int_{\GammaEndoR} \PRV \,J \mecF^{-T}\mecNact \, d\Gamma_0
	\qquad \text{on $\GammaBaseR$} \\
	\end{split}
	\right.
	\end{equation}

	\item \textbf{Weighted stress distribution}.
	Finally, we consider the case in which we set $\BCweight = \XiHat$ (as defined in Sec.~\ref{subsec: fibers}).
	The function $\XiHat$ is defined such that we have $\XiHat \simeq 1$ on $\GammaBaseL$, $\XiHat \simeq 0$ on $\GammaBaseR$ and we have a smooth transition on the septum.
	With this choice, the energy-consistent boundary condition of Eq.~\eqref{eqn:BC_generic} reads
	\begin{equation}
	\begin{split}
	\tenspiola\mecNact
	=
	| J \mecF^{-T} \mecNref | \left[ \PLV(t)\BCmecVbaseLV(t,\XiHat) + \PRV(t)\BCmecVbaseRV(t,\XiHat) \right] ,
	\end{split}
	\end{equation}
	having defined the vectors $\BCmecVbaseLV$ and $\BCmecVbaseRV$ as in Eq.~\eqref{eqn: weighted_bc}.
\end{itemize}

Based upon our experience, the uniform stress distribution approach does not typically provide meaningful results. Indeed, since the stress is redistributed on the whole base without accounting for the closeness to the two chambers, a net angular momentum results on the elastic body, making it rotate during systole.
Conversely, both the uniform stress distribution approach over each base and the weighted stress distribution approach overcome this issue, thanks to a more realistic distribution of the stress.
While the two strategies globally provide very similar results, the latter allows for a smoother solution close to the interface between the left and right bases.
For this reason, in this paper we focus on the weighted stress distribution approach.

%% file: parts_app_schur.tex
\section{3D-0D saddle-point problem resolution}
\label{app:schur}
We solve the non-linear saddle-point problem~\eqref{eqn: saddle-point} by means of the following Newton algorithm (where the subscript $n+1$ is understood):
\begin{itemize}
	\item We set $\displdofs_{\text{h}}^{(0)} = \displdofs_{\text{h}}^{n}$,  $\PLV^{(0)} = \PLV^{n}$ and $\PRV^{(0)} = \PRV^{n}$
	\item For $j=1,2,\dots$, until convergence, we solve the linear system
	\begin{equation}\label{eqn: Msi newton:biventricle}
		\begin{pmatrix}
			J_{\displ, \displ}^{(j-1)} & J_{\displ, \PLV}^{(j-1)} & J_{\displ, \PRV}^{(j-1)} \\
			J_{\PLV, \displ}^{(j-1)} & 0 & 0\\
			J_{\PRV, \displ}^{(j-1)} & 0 & 0\\
		\end{pmatrix}
		\begin{pmatrix}
			\Delta{\displdofs_{\text{h}}^{(j)}} \\
			\Delta{\PLV^{(j)}} \\
			\Delta{\PRV^{(j)}} \\
		\end{pmatrix}
		=
		\begin{pmatrix}
			\mathbf{r}_{\displ}^{(j-1)} \\
			r_{\PLV}^{(j-1)} \\
			r_{\PRV}^{(j-1)} \\
		\end{pmatrix},
	\end{equation}
	where
	
	$J_{\displ, \displ}^{(j-1)} = \frac{\partial}{\partial\displdofs} \mathbf{r}_{\displ}(\displdofs_{\text{h}}^{(j-1)},\PLV^{(j-1)},\PRV^{(j-1)})$,
	
	$J_{\displ, \PLV}^{(j-1)} = \frac{\partial}{\partial \PLV} \mathbf{r}_{\displ}(\displdofs_{\text{h}}^{(j-1)},\PLV^{(j-1)},\PRV^{(j-1)}), \quad
	J_{\displ, \PRV}^{(j-1)} = \frac{\partial}{\partial \PRV} \mathbf{r}_{\displ}(\displdofs_{\text{h}}^{(j-1)},\PLV^{(j-1)},\PRV^{(j-1)})$,
	
	$J_{\PLV, \displ}^{(j-1)} = \frac{\partial}{\partial\displdofs} r_{\PLV}(\displdofs_{\text{h}}^{(j-1)}), \quad
	J_{\PRV, \displ}^{(j-1)} = \frac{\partial}{\partial\displdofs} r_{\PRV}(\displdofs_{\text{h}}^{(j-1)})$,
	
	
	\item We update
	
	$\displdofs_{\text{h}}^{(j)} = \displdofs_{\text{h}}^{(j-1)} + \Delta{\displdofs_{\text{h}}^{(j)}}$,
	$\PLV^{(j)} = \PLV^{(j-1)} + \Delta{\PLV^{(j)}}$ and
	$\PRV^{(j)} = \PRV^{(j-1)} + \Delta{\PRV^{(j)}}.$
	\item When the convergence criterion (based on the increment) is satisfied, we set
	
	$\displdofs_{\text{h}}^{n+1} = \displdofs_{\text{h}}^{(j)}$, $\PLV^{n+1} = \PLV^{(j)}$ and $\PRV^{n+1} = \PRV^{(j)}.$
\end{itemize}

\noindent
We solve the saddle-point problem~\eqref{eqn: Msi newton:biventricle} via Schur complement reduction \cite{benzi2005numerical}. Specifically, system~\eqref{eqn: saddle-point} can be written as
\begin{equation}\label{eqn: saddle-point_system}
	\begin{cases}
		&J_{\displ, \displ} \Delta{\displdofs_{\text{h}}} + J_{\PLV, \displ} \Delta{\PLV} + J_{\PRV, \displ} \Delta{\PRV} = \mathbf{r}_{\displ}  \\
		& J_{\PLV, \displ} \Delta{\displdofs_{\text{h}}} = r_{\PLV}  \\
		& J_{\PRV, \displ} \Delta{\displdofs_{\text{h}}} = r_{\PRV}
	\end{cases}
\end{equation}
where for simplicity we omit the superscript $(j)$. Deriving $\Delta{\displdofs_{\text{h}}}$ form the first equation of~\eqref{eqn: saddle-point_system} we have
\begin{equation}
	\begin{cases}\label{eqn: saddle-point_system2}
		&\Delta{\displdofs_{\text{h}}} = \mathbf{v} - \mathbf{w}_{\text{L}} \Delta{\PLV} - \mathbf{w}_{\text{R}} \Delta{\PRV} \\
		& \alpha_{\text{LL}} \Delta{\PLV} +
		\alpha_{\text{LR}}\Delta{\PRV} = b_{\text{L}} \\
		& \alpha_{\text{RL}} \Delta{\PLV} +
		\alpha_{\text{RR}} \Delta{\PRV} = b_{\text{R}}
		\\
	\end{cases}
\end{equation}
where
$$ \alpha_{\text{LL}} = J_{\PLV, \displ} \mathbf{w}_{\text{L}}, \qquad
\alpha_{\text{LR}} = J_{\PLV, \displ} \mathbf{w}_{\text{R}}, \qquad
\alpha_{\text{RL}} = J_{\PRV, \displ} \mathbf{w}_{\text{L}}, \qquad
\alpha_{\text{RR}} = J_{\PRV, \displ} \mathbf{w}_{\text{R}} ,$$
$$ b_{\text{L}}  = J_{\PLV, \displ} \mathbf{v}  - r_{\PLV}, \qquad b_{\text{R}} = J_{\PRV, \displ} \mathbf{v}
- r_{\PRV}, $$
with
\begin{equation}\label{eqn: three_systems}
	\mathbf{w}_{\text{L}} = J_{\displ, \displ}^{-1} J_{\PLV, \displ} \qquad \mathbf{w}_{\text{R}} = J_{\displ, \displ}^{-1}J_{\PRV, \displ}, \qquad \mathbf{v} = J_{\displ, \displ}^{-1}\mathbf{r}_{\displ}.
\end{equation}
Solving equation~\eqref{eqn: saddle-point_system2} we obtain
\begin{equation}\label{eqn: schur_solution_2}	
\Delta{\displdofs_{\text{h}}} = \mathbf{v} - \mathbf{w}_{\text{L}} \Delta{\PLV} - \mathbf{w}_{\text{R}} \Delta{\PRV},$$ 
$$\Delta{\PLV} = \frac{b_{\text{L}}\alpha_{\text{RR}} + b_{\text{R}} \alpha_{\text{LR}}}{\alpha_{\text{LL}}\alpha_{\text{RR}} - \alpha_{\text{RL}} \alpha_{\text{LR}}}, \quad
\Delta{\PRV} = \frac{b_{\text{R}}\alpha_{\text{LL}} + b_{\text{L}} \alpha_{\text{RL}}}{\alpha_{\text{LL}}\alpha_{\text{RR}} - \alpha_{\text{RL}} \alpha_{\text{LR}}}.
\end{equation}
Notice that we have to the solve three linear systems~\eqref{eqn: three_systems} in order to obtain the solution (20).

%% file: EMBiventricular_CMAME.bbl
\begin{thebibliography}{100}
\expandafter\ifx\csname url\endcsname\relax
  \def\url#1{\texttt{#1}}\fi
\expandafter\ifx\csname urlprefix\endcsname\relax\def\urlprefix{URL }\fi
\expandafter\ifx\csname href\endcsname\relax
  \def\href#1#2{#2} \def\path#1{#1}\fi

\bibitem{gurev2011models}
V.~Gurev, T.~Lee, J.~Constantino, H.~Arevalo, N.~Trayanova, Models of cardiac
  electromechanics based on individual hearts imaging data, Biomechanics and
  Modeling in Mechanobiology 10~(3) (2011) 295--306.

\bibitem{augustin2016anatomically}
C.~Augustin, A.~Neic, M.~Liebmann, A.~Prassl, S.~Niederer, G.~Haase, G.~Plank,
  Anatomically accurate high resolution modeling of human whole heart
  electromechanics: a strongly scalable algebraic multigrid solver method for
  nonlinear deformation, Journal of Computational Physics 305 (2016) 622--646.

\bibitem{augustin2020impact}
C.~Augustin, T.~Fastl, A.~Neic, C.~Bellini, J.~Whitaker, R.~Rajani,
  M.~O’Neill, M.~Bishop, G.~Plank, S.~Niederer, The impact of wall thickness
  and curvature on wall stress in patient-specific electromechanical models of
  the left atrium, Biomechanics and Modeling in Mechanobiology 19~(3) (2020)
  1015--1034.

\bibitem{land2018influence}
S.~Land, S.~Niederer, Influence of atrial contraction dynamics on cardiac
  function, International Journal for Numerical Methods in Biomedical
  Engineering 34~(3) (2018) e2931.

\bibitem{nordsletten2011coupling}
D.~Nordsletten, S.~Niederer, M.~Nash, P.~Hunter, N.~Smith, Coupling
  multi-physics models to cardiac mechanics, Progress in Biophysics and
  Molecular Biology 104 (2011) 77--88.

\bibitem{strocchi2020publicly}
M.~Strocchi, C.~Augustin, M.~Gsell, E.~Karabelas, A.~Neic, K.~Gillette,
  O.~Razeghi, A.~Prassl, E.~Vigmond, J.~Behar, et~al., A publicly available
  virtual cohort of four-chamber heart meshes for cardiac electro-mechanics
  simulations, PloS One 15 (2020) e0235145.

\bibitem{math9111247}
T.~Gerach, S.~Schuler, J.~Fröhlich, L.~Lindner, E.~Kovacheva, R.~Moss,
  E.~Wülfers, G.~Seemann, C.~Wieners, A.~Loewe, Electro-mechanical whole-heart
  digital twins: A fully coupled multi-physics approach, Mathematics 9~(11)
  (2021).

\bibitem{sermesant2012patient}
M.~Sermesant, R.~Chabiniok, P.~Chinchapatnam, T.~Mansi, F.~Billet, P.~Moireau,
  J.~Peyrat, K.~Wong, J.~Relan, K.~Rhode, et~al., Patient-specific
  electromechanical models of the heart for the prediction of pacing acute
  effects in crt: a preliminary clinical validation, Medical Image Analysis
  16~(1) (2012) 201--215.

\bibitem{peirlinck2021precision}
M.~Peirlinck, F.~Costabal, J.~Yao, J.~Guccione, S.~Tripathy, Y.~Wang,
  D.~Ozturk, P.~Segars, T.~Morrison, S.~Levine, et~al., Precision medicine in
  human heart modeling, Biomechanics and Modeling in Mechanobiology (2021)
  1--29.

\bibitem{smith2004multiscale}
N.~Smith, D.~Nickerson, E.~Crampin, P.~Hunter, Multiscale computational
  modelling of the heart, Acta Numerica 13 (2004) 371.

\bibitem{chabiniok2016multiphysics}
R.~Chabiniok, V.~Wang, M.~Hadjicharalambous, L.~Asner, J.~Lee, M.~Sermesant,
  E.~Kuhl, A.~Young, P.~Moireau, M.~Nash, et~al., Multiphysics and multiscale
  modelling, data--model fusion and integration of organ physiology in the
  clinic: ventricular cardiac mechanics, Interface Focus 6~(2) (2016) 20150083.

\bibitem{crampin2004computational}
E.~Crampin, M.~Halstead, P.~Hunter, P.~Nielsen, D.~Noble, N.~Smith, M.~Tawhai,
  Computational physiology and the physiome project, Experimental Physiology
  89~(1) (2004) 1--26.

\bibitem{marx2020personalization}
L.~Marx, M.~Gsell, A.~Rund, F.~Caforio, A.~Prassl, G.~Toth-Gayor, T.~Kuehne,
  C.~Augustin, G.~Plank, Personalization of electro-mechanical models of the
  pressure-overloaded left ventricle: fitting of windkessel-type afterload
  models, Philosophical Transactions of the Royal Society A 378~(2173) (2020)
  20190342.

\bibitem{Gerbi}
A.~Gerbi, L.~Dede', A.~Quarteroni, A monolithic algorithm for the simulation of
  cardiac electromechanics in the human left ventricle, Mathematics in
  Engineering 1 (2018).

\bibitem{regazzoni2020model}
F.~Regazzoni, M.~Salvador, P.~C. Africa, M.~Fedele, L.~Dede', A.~Quarteroni, A
  cardiac electromechanics model coupled with a lumped parameters model for
  closed-loop blood circulation. part i: model derivation (2020).
\newblock \href {http://arxiv.org/abs/2011.15040} {\path{arXiv:2011.15040}}.

\bibitem{salvador2020intergrid}
M.~Salvador, L.~Ded{\`e}, A.~Quarteroni, An intergrid transfer operator using
  radial basis functions with application to cardiac electromechanics,
  Computational Mechanics 66 (2020) 491--511.

\bibitem{levrero2020sensitivity}
F.~Levrero-Florencio, F.~Margara, E.~Zacur, A.~Bueno-Orovio, Z.~Wang,
  A.~Santiago, J.~Aguado-Sierra, G.~Houzeaux, V.~Grau, D.~Kay, et~al.,
  Sensitivity analysis of a strongly-coupled human-based electromechanical
  cardiac model: Effect of mechanical parameters on physiologically relevant
  biomarkers, Computer Methods in Applied Mechanics and Engineering 361 (2020)
  112762.

\bibitem{propp2020orthotropic}
A.~Propp, A.~Gizzi, F.~Levrero-Florencio, R.~Ruiz-Baier, An orthotropic
  electro-viscoelastic model for the heart with stress-assisted diffusion,
  Biomechanics and Modeling in Mechanobiology 19~(2) (2020) 633--659.

\bibitem{rossi2014thermodynamically}
S.~Rossi, T.~Lassila, R.~Ruiz-Baier, A.~Sequeira, A.~Quarteroni,
  Thermodynamically consistent orthotropic activation model capturing
  ventricular systolic wall thickening in cardiac electromechanics, European
  Journal of Mechanics-A/Solids 48 (2014) 129--142.

\bibitem{palit2015computational}
A.~Palit, S.~Bhudia, T.~Arvanitis, G.~Turley, M.~Williams, Computational
  modelling of left-ventricular diastolic mechanics: Effect of fibre
  orientation and right-ventricle topology, Journal of Biomechanics 48~(4)
  (2015) 604--612.

\bibitem{sermesant2005simulation}
M.~Sermesant, K.~Rhode, G.~Sanchez-Ortiz, O.~Camara, R.~Andriantsimiavona,
  S.~Hegde, D.~Rueckert, P.~Lambiase, C.~Bucknall, E.~Rosenthal, et~al.,
  Simulation of cardiac pathologies using an electromechanical biventricular
  model and xmr interventional imaging, Medical Image Analysis 9~(5) (2005)
  467--480.

\bibitem{chapelle2009numerical}
D.~Chapelle, M.~Fern{\'a}ndez, J.~Gerbeau, P.~Moireau, J.~Sainte-Marie,
  N.~Zemzemi, Numerical simulation of the electromechanical activity of the
  heart, in: International Conference on Functional Imaging and Modeling of the
  Heart, Springer, 2009, pp. 357--365.

\bibitem{goktepe2010electromechanics}
S.~G{\"o}ktepe, E.~Kuhl, Electromechanics of the heart: a unified approach to
  the strongly coupled excitation--contraction problem, Computational Mechanics
  45~(2) (2010) 227--243.

\bibitem{crozier2016image}
A.~Crozier, C.~Augustin, A.~Neic, A.~Prassl, M.~Holler, T.~Fastl, A.~Hennemuth,
  K.~Bredies, T.~Kuehne, M.~Bishop, et~al., Image-based personalization of
  cardiac anatomy for coupled electromechanical modeling, Annals of Biomedical
  Engineering 44~(1) (2016) 58--70.

\bibitem{hirschvogel2017monolithic}
M.~Hirschvogel, M.~Bassilious, L.~Jagschies, S.~Wildhirt, M.~Gee, A monolithic
  3d-0d coupled closed-loop model of the heart and the vascular system:
  experiment-based parameter estimation for patient-specific cardiac mechanics,
  International Journal for Numerical Methods in Biomedical Engineering 33
  (2017) e2842.

\bibitem{ahmad2018multiphysics}
A.~Ahmad~Bakir, A.~Al~Abed, M.~Stevens, N.~Lovell, S.~Dokos, A multiphysics
  biventricular cardiac model: Simulations with a left-ventricular assist
  device, Frontiers in Physiology 9 (2018) 1259.

\bibitem{garcia2019towards}
E.~Garcia-Blanco, R.~Ortigosa, A.~Gil, J.~Bonet, Towards an efficient
  computational strategy for electro-activation in cardiac mechanics, Computer
  Methods in Applied Mechanics and Engineering 356 (2019) 220--260.

\bibitem{augustin2020physiologically}
C.~Augustin, M.~Gsell, E.~Karabelas, G.~Plank, Physiologically valid 3d-0d
  closed loop model of the heart and circulation--modeling the acute response
  to altered loading and contractility, arXiv preprint arXiv:2009.08802 (2020).

\bibitem{kerckhoffs2007coupling}
R.~Kerckhoffs, M.~Neal, Q.~Gu, J.~Bassingthwaighte, J.~Omens, A.~McCulloch,
  Coupling of a 3d finite element model of cardiac ventricular mechanics to
  lumped systems models of the systemic and pulmonic circulation, Annals of
  Biomedical Engineering 35 (2007) 1--18.

\bibitem{dede2021modeling}
L.~Ded{\`e}, F.~Regazzoni, C.~Vergara, P.~Zunino, M.~Guglielmo, R.~Scrofani,
  L.~Fusini, C.~Cogliati, G.~Pontone, A.~Quarteroni, Modeling the cardiac
  response to hemodynamic changes associated with {COVID}-19: a computational
  study, Mathematical Biosciences and Engineering 18~(4) (2021) 3364--3383.

\bibitem{guan2020effect}
D.~Guan, J.~Yao, X.~Luo, H.~Gao, Effect of myofibre architecture on ventricular
  pump function by using a neonatal porcine heart model: from dt-mri to
  rule-based methods, Royal Society Open Science 7~(4) (2020) 191655.

\bibitem{wang2021human}
Z.~Wang, A.~Santiago, X.~Zhou, L.~Wang, F.~Margara, F.~Levrero-Florencio,
  A.~Das, C.~Kelly, E.~Dall'Armellina, M.~Vazquez, et~al., Human biventricular
  electromechanical simulations on the progression of electrocardiographic and
  mechanical abnormalities in post-myocardial infarction, EP Europace
  23~(Supplement\_1) (2021) i143--i152.

\bibitem{elzinga1973pressure}
G.~Elzinga, N.~Westerhof, Pressure and flow generated by the left ventricle
  against different impedances, Circulation Research 32~(2) (1973) 178--186.

\bibitem{liu2015multi}
H.~Liu, F.~Liang, J.~Wong, T.~Fujiwara, W.~Ye, K.~Tsubota, M.~Sugawara,
  Multi-scale modeling of hemodynamics in the cardiovascular system, Acta
  Mechanica Sinica 31~(4) (2015) 446--464.

\bibitem{segers2008three}
P.~Segers, E.~Rietzschel, M.~De~Buyzere, N.~Stergiopulos, N.~Westerhof,
  L.~Van~Bortel, T.~Gillebert, P.~Verdonck, Three-and four-element windkessel
  models: assessment of their fitting performance in a large cohort of healthy
  middle-aged individuals, Proceedings of the Institution of Mechanical
  Engineers, Part H: Journal of Engineering in Medicine 222~(4) (2008)
  417--428.

\bibitem{stergiopulos1999total}
N.~Stergiopulos, B.~Westerhof, N.~Westerhof, Total arterial inertance as the
  fourth element of the {W}indkessel model, American Journal of
  Physiology-Heart and Circulatory Physiology 276~(1) (1999) H81--H88.

\bibitem{wang2003time}
J.~Wang, A.~O'Brien, N.~Shrive, K.~Parker, J.~Tyberg, Time-domain
  representation of ventricular-arterial coupling as a windkessel and wave
  system, American Journal of Physiology-Heart and Circulatory Physiology
  284~(4) (2003) H1358--H1368.

\bibitem{westerhof1991normalized}
N.~Westerhof, G.~Elzinga, Normalized input impedance and arterial decay time
  over heart period are independent of animal size, American Journal of
  Physiology-Regulatory, Integrative and Comparative Physiology 261~(1) (1991)
  R126--R133.

\bibitem{dede2020segregated}
L.~Ded{\`e}, A.~Gerbi, A.~Quarteroni, Segregated algorithms for the numerical
  simulation of cardiac electromechanics in the left human ventricle, in:
  D.~Ambrosi, P.~Ciarletta (Eds.), The Mathematics of Mechanobiology, Springer,
  2020, pp. 81--116.

\bibitem{eriksson2013influence}
T.~Eriksson, A.~J. Prassl, G.~Plank, G.~Holzapfel, Influence of myocardial
  fiber/sheet orientations on left ventricular mechanical contraction,
  Mathematics and Mechanics of Solids 18~(6) (2013) 592--606.

\bibitem{usyk2002computational}
T.~Usyk, I.~LeGrice, A.~McCulloch, Computational model of three-dimensional
  cardiac electromechanics, Computing and Visualization in Science 4~(4) (2002)
  249--257.

\bibitem{blanco20103d}
P.~Blanco, R.~Feij{\'o}o, et~al., A 3d-1d-0d computational model for the entire
  cardiovascular system, Computational Mechanics 29 (2010) 5887--5911.

\bibitem{arts2005adaptation}
T.~Arts, T.~Delhaas, P.~Bovendeerd, X.~Verbeek, F.~Prinzen, Adaptation to
  mechanical load determines shape and properties of heart and circulation: the
  circadapt model, American Journal of Physiology-Heart and Circulatory
  Physiology 288 (2005) H1943--H1954.

\bibitem{neal2007subject}
M.~Neal, J.~Bassingthwaighte, Subject-specific model estimation of cardiac
  output and blood volume during hemorrhage, Cardiovascular Engineering 7~(3)
  (2007) 97--120.

\bibitem{paeme2011mathematical}
S.~Paeme, K.~Moorhead, J.~Chase, B.~Lambermont, P.~Kolh, V.~D'orio, L.~Pierard,
  M.~Moonen, P.~Lancellotti, P.~Dauby, et~al., Mathematical multi-scale model
  of the cardiovascular system including mitral valve dynamics. application to
  ischemic mitral insufficiency, Biomedical {E}ngineering {O}nline 10~(1)
  (2011) 1--20.

\bibitem{regazzoni2020numerical}
F.~Regazzoni, M.~Salvador, P.~C. Africa, M.~Fedele, L.~Dede', A.~Quarteroni, A
  cardiac electromechanics model coupled with a lumped parameters model for
  closed-loop blood circulation. part ii: numerical approximation (2020).
\newblock \href {http://arxiv.org/abs/2011.15051} {\path{arXiv:2011.15051}}.

\bibitem{piersanti2021modeling}
R.~Piersanti, P.~Africa, M.~Fedele, C.~Vergara, L.~Ded{\`e}, A.~Corno,
  A.~Quarteroni, Modeling cardiac muscle fibers in ventricular and atrial
  electrophysiology simulations, Computer Methods in Applied Mechanics and
  Engineering 373 (2021) 113468.

\bibitem{streeter1969fiber}
D.~D. Streeter~Jr, H.~M. Spotnitz, D.~P. Patel, J.~Ross~Jr, E.~H. Sonnenblick,
  Fiber orientation in the canine left ventricle during diastole and systole,
  Circulation Research 24~(3) (1969) 339--347.

\bibitem{roberts1979influence}
D.~E. Roberts, L.~T. Hersh, A.~M. Scher, Influence of cardiac fiber orientation
  on wavefront voltage, conduction velocity, and tissue resistivity in the
  dog., Circulation Research 44~(5) (1979) 701--712.

\bibitem{gil2019influence}
D.~Gil, R.~Aris, A.~Borras, E.~Ram{\'\i}rez, R.~Sebastian, M.~V{\'a}zquez,
  Influence of fiber connectivity in simulations of cardiac biomechanics,
  International Journal of Computer Assisted Radiology and Surgery 14~(1)
  (2019) 63--72.

\bibitem{carreras2012left}
F.~Carreras, J.~Garcia-Barnes, D.~Gil, S.~Pujadas, C.~Li, R.~Suarez-Arias,
  R.~Leta, X.~Alomar, M.~Ballester, G.~Pons-Llado, Left ventricular torsion and
  longitudinal shortening: two fundamental components of myocardial mechanics
  assessed by tagged cine-mri in normal subjects, The {I}nternational {J}ournal
  of {C}ardiovascular {I}maging 28~(2) (2012) 273--284.

\bibitem{bayer2012novel}
J.~Bayer, R.~Blake, G.~Plank, N.~Trayanova, A novel rule-based algorithm for
  assigning myocardial fiber orientation to computational heart models, Annals
  of Biomedical Engineering 40~(10) (2012) 2243--2254.

\bibitem{wong2014generating}
J.~Wong, E.~Kuhl, {Generating fibre orientation maps in human heart models
  using Poisson interpolation}, Computer Methods in Biomechanics and Biomedical
  Engineering 17~(11) (2014) 1217--1226.

\bibitem{doste2019rule}
R.~Doste, D.~Soto-Iglesias, G.~Bernardino, A.~Alcaine, R.~Sebastian,
  S.~Giffard-Roisin, M.~Sermesant, A.~Berruezo, D.~Sanchez-Quintana, O.~Camara,
  A rule-based method to model myocardial fiber orientation in cardiac
  biventricular geometries with outflow tracts, International Journal for
  Numerical Methods in Biomedical Engineering 35~(4) (2019) e3185.

\bibitem{quarteroni2017integrated}
A.~Quarteroni, T.~Lassila, S.~Rossi, R.~Ruiz-Baier, {Integrated
  Heart—Coupling multiscale and multiphysics models for the simulation of the
  cardiac function}, Computer Methods in Applied Mechanics and Engineering 314
  (2017) 345--407.

\bibitem{azzolin2020effect}
L.~Azzolin, L.~Ded{\`e}, A.~Gerbi, A.~Quarteroni, Effect of fibre orientation
  and bulk modulus on the electromechanical modelling of human ventricles,
  Mathematics in Engineering 2~(4) (2020) 614--638.

\bibitem{pluijmert2017determinants}
M.~Pluijmert, T.~Delhaas, A.~De~la Parra, W.~Kroon, F.~Prinzen, P.~Bovendeerd,
  Determinants of biventricular cardiac function: a mathematical model study on
  geometry and myofiber orientation, Biomechanics and Modeling in
  Mechanobiology 16~(2) (2017) 721--729.

\bibitem{guan2021modelling}
D.~Guan, X.~Zhuan, W.~Holmes, X.~Luo, H.~Gao, Modelling of fibre dispersion and
  its effects on cardiac mechanics from diastole to systole, Journal of
  Engineering Mathematics 128~(1) (2021) 1--24.

\bibitem{ahmad2018region}
F.~Ahmad, S.~Soe, N.~White, R.~Johnston, I.~Khan, J.~Liao, M.~Jones, R.~Prabhu,
  I.~Maconochie, P.~Theobald, Region-specific microstructure in the neonatal
  ventricles of a porcine model, Annals of Biomedical Engineering 46~(12)
  (2018) 2162--2176.

\bibitem{sommer2015biomechanical}
G.~Sommer, A.~Schriefl, M.~Andr{\"a}, M.~Sacherer, C.~Viertler, H.~Wolinski,
  G.~Holzapfel, Biomechanical properties and microstructure of human
  ventricular myocardium, Acta Biomaterialia 24 (2015) 172--192.

\bibitem{lin1998multiaxial}
D.~Lin, F.~Yin, A multiaxial constitutive law for mammalian left ventricular
  myocardium in steady-state barium contracture or tetanus, Journal of
  Biomechanical Engineering 120~(4) (1998) 504--517.

\bibitem{genet2014distribution}
M.~Genet, L.~Lee, R.~Nguyen, H.~Haraldsson, G.~Acevedo-Bolton, Z.~Zhang, L.~Ge,
  K.~Ordovas, S.~Kozerke, J.~Guccione, Distribution of normal human left
  ventricular myofiber stress at end diastole and end systole: a target for in
  silico design of heart failure treatments, Journal of Applied Physiology
  117~(2) (2014) 142--152.

\bibitem{sack2018construction}
K.~Sack, E.~Aliotta, D.~Ennis, J.~Choy, G.~Kassab, J.~Guccione, T.~Franz,
  Construction and validation of subject-specific biventricular finite-element
  models of healthy and failing swine hearts from high-resolution dt-mri,
  Frontiers in {P}hysiology 9 (2018) 539.

\bibitem{wenk2012first}
J.~Wenk, D.~Klepach, L.~Lee, Z.~Zhang, L.~Ge, E.~Tseng, A.~Martin, S.~Kozerke,
  J.~Gorman~III, R.~Gorman, et~al., First evidence of depressed contractility
  in the border zone of a human myocardial infarction, The Annals of Thoracic
  Surgery 93~(4) (2012) 1188--1193.

\bibitem{eriksson2013modeling}
T.~Eriksson, A.~Prassl, G.~Plank, G.~Holzapfel, Modeling the dispersion in
  electromechanically coupled myocardium, International Journal for Numerical
  Methods in Biomedical Engineering 29~(11) (2013) 1267--1284.

\bibitem{maceira2006normalized}
A.~Maceira, S.~Prasad, M.~Khan, D.~Pennell, Normalized left ventricular
  systolic and diastolic function by steady state free precession
  cardiovascular magnetic resonance, Journal of Cardiovascular Magnetic
  Resonance 8~(3) (2006) 417--426.

\bibitem{tamborini2010reference}
G.~Tamborini, N.~Marsan, P.~Gripari, F.~Maffessanti, D.~Brusoni, M.~Muratori,
  E.~Caiani, C.~Fiorentini, M.~Pepi, Reference values for right ventricular
  volumes and ejection fraction with real-time three-dimensional
  echocardiography: evaluation in a large series of normal subjects, Journal of
  the American Society of Echocardiography 23~(2) (2010) 109--115.

\bibitem{maceira2006reference}
A.~Maceira, S.~Prasad, M.~Khan, D.~Pennell, Reference right ventricular
  systolic and diastolic function normalized to age, gender and body surface
  area from steady-state free precession cardiovascular magnetic resonance,
  European Heart Journal 27~(23) (2006) 2879--2888.

\bibitem{sugimoto2017echocardiographic}
T.~Sugimoto, R.~Dulgheru, A.~Bernard, F.~Ilardi, L.~Contu, K.~Addetia,
  L.~Caballero, N.~Akhaladze, G.~Athanassopoulos, D.~Barone, et~al.,
  Echocardiographic reference ranges for normal left ventricular 2d strain:
  results from the eacvi norre study, European Heart Journal-Cardiovascular
  Imaging 18~(8) (2017) 833--840.

\bibitem{bishop1997clinical}
A.~Bishop, P.~White, P.~Oldershaw, R.~Chaturvedi, C.~Brookes, A.~Redington,
  Clinical application of the conductance catheter technique in the adult human
  right ventricle, International Journal of Cardiology 58~(3) (1997) 211--221.

\bibitem{emilsson2006mitral}
K.~Emilsson, R.~Egerlid, B.~Nygren, B.~Wandt, Mitral annulus motion versus
  long-axis fractional shortening, Experimental \& Clinical Cardiology 11~(4)
  (2006) 302.

\bibitem{sechtem1987regional}
U.~Sechtem, B.~Sommerhoff, W.~Markiewicz, R.~White, M.~Cheitlin, C.~Higgins,
  Regional left ventricular wall thickening by magnetic resonance imaging:
  evaluation in normal persons and patients with global and regional
  dysfunction, The American Journal of Cardiology 59~(1) (1987) 145--151.

\bibitem{lee2019rule}
A.~Lee, U.~Nguyen, O.~Razeghi, J.~Gould, B.~Sidhu, B.~Sieniewicz, J.~Behar,
  M.~Mafi-Rad, G.~Plank, F.~Prinzen, et~al., A rule-based method for predicting
  the electrical activation of the heart with cardiac resynchronization therapy
  from non-invasive clinical data, Medical {I}mage {A}nalysis 57 (2019)
  197--213.

\bibitem{bayer2005laplace}
J.~Bayer, J.~Beaumont, A.~Krol, {Laplace--Dirichlet energy field specification
  for deformable models. An FEM approach to active contour fitting}, Annals of
  Biomedical Engineering 33~(9) (2005) 1175--1186.

\bibitem{franzone2014mathematical}
P.~Franzone, L.~Pavarino, S.~Scacchi, Mathematical Cardiac Electrophysiology,
  Vol.~13, Springer, 2014.

\bibitem{colli2018numerical}
P.~Franzone, L.~Pavarino, S.~Scacchi, A numerical study of scalable cardiac
  electro-mechanical solvers on {HPC} architectures, Frontiers in {P}hysiology
  9 (2018) 268.

\bibitem{nobile2012active}
F.~Nobile, A.~Quarteroni, R.~Ruiz-Baier, An active strain electromechanical
  model for cardiac tissue, International Journal for Numerical Methods in
  Biomedical Engineering 28 (2012) 52--71.

\bibitem{niederer2011length}
S.~Niederer, G.~Plank, P.~Chinchapatnam, M.~Ginks, P.~Lamata, K.~Rhode,
  C.~Rinaldi, R.~Razavi, N.~Smith, Length-dependent tension in the failing
  heart and the efficacy of cardiac resynchronization therapy, Cardiovascular
  Research 89 (2011) 336--343.

\bibitem{ruiz2014mathematical}
R.~Ruiz-Baier, A.~Gizzi, S.~Rossi, C.~Cherubini, A.~Laadhari, S.~Filippi,
  A.~Quarteroni, Mathematical modelling of active contraction in isolated
  cardiomyocytes, Mathematical Medicine and Biology: a Journal of the IMA 31
  (2014) 259--283.

\bibitem{land2017model}
S.~Land, S.~Park-Holohan, N.~Smith, C.~Dos~Remedios, J.~Kentish, S.~Niederer, A
  model of cardiac contraction based on novel measurements of tension
  development in human cardiomyocytes, Journal of {M}olecular and {C}ellular
  {C}ardiology 106 (2017) 68--83.

\bibitem{regazzoni2020biophysically}
F.~Regazzoni, L.~Ded{\`e}, A.~Quarteroni, Biophysically detailed mathematical
  models of multiscale cardiac active mechanics, PLoS {C}omputational {B}iology
  16 (2020) e1008294.

\bibitem{regazzoni2020machine}
F.~Regazzoni, L.~Ded{\`e}, A.~Quarteroni, Machine learning of multiscale active
  force generation models for the efficient simulation of cardiac
  electromechanics, Computer Methods in Applied Mechanics and Engineering 370
  (2020) 113268.

\bibitem{guccione1991passive}
J.~Guccione, A.~McCulloch, L.~Waldman, Passive material properties of intact
  ventricular myocardium determined from a cylindrical model, Journal of
  Biomechanical Engineering 113 (1991) 42--55.

\bibitem{guccione1991finite}
J.~Guccione, A.~McCulloch, Finite element modeling of ventricular mechanics,
  in: Theory of Heart, Springer, 1991, pp. 121--144.

\bibitem{ogden1997non}
R.~Ogden, Non-linear elastic deformations, Courier Corporation, 1997.

\bibitem{holzapfel2009constitutive}
G.~Holzapfel, R.~Ogden, Constitutive modelling of passive myocardium: a
  structurally based framework for material characterization, Philosophical
  Transactions of the Royal Society A: Mathematical, Physical and Engineering
  Sciences 367 (2009) 3445--3475.

\bibitem{quarteroni2016geometric}
A.~Quarteroni, A.~Veneziani, C.~Vergara, Geometric multiscale modeling of the
  cardiovascular system, between theory and practice, Computer Methods in
  Applied Mechanics and Engineering 302 (2016) 193--252.

\bibitem{vergara2014patient}
C.~Vergara, S.~Palamara, D.~Catanzariti, F.~Nobile, E.~Faggiano, C.~Pangrazzi,
  M.~Centonze, M.~Maines, A.~Quarteroni, G.~Vergara, Patient-specific
  generation of the {P}urkinje network driven by clinical measurements of a
  normal propagation, Medical \& {B}iological {E}ngineering \& {C}omputing 52
  (2014) 813--826.

\bibitem{costabal2016generating}
F.~Costabal, D.~Hurtado, E.~Kuhl, Generating purkinje networks in the human
  heart, Journal of Biomechanics 49 (2016) 2455--2465.

\bibitem{landajuela2018numerical}
M.~Landajuela, C.~Vergara, A.~Gerbi, L.~Ded{\`e}, L.~Formaggia, A.~Quarteroni,
  Numerical approximation of the electromechanical coupling in the left
  ventricle with inclusion of the purkinje network, International {J}ournal for
  {N}umerical {M}ethods in {B}iomedical {E}ngineering 34 (2018) e2984.

\bibitem{ten2006alternans}
K.~ten Tusscher, A.~Panfilov, Alternans and spiral breakup in a human
  ventricular tissue model, American Journal of Physiology-Heart and
  Circulatory Physiology 291 (2006) H1088--H1100.

\bibitem{bers2001excitation}
D.~Bers, Excitation-contraction coupling and cardiac contractile force, Vol.
  237, Springer Science \& Business Media, 2001.

\bibitem{peng1997compressible}
S.~Peng, W.~Chang, A compressible approach in finite element analysis of
  rubber-elastic materials, Computers \& Structures 62 (1997) 573--593.

\bibitem{doll2000development}
S.~Doll, K.~Schweizerhof, On the development of volumetric strain energy
  functions, Journal of Applied Mechanics 67 (2000) 17--21.

\bibitem{Gerbi2018monolithic}
A.~Gerbi, L.~Ded{\`e}, A.~Quarteroni, A monolithic algorithm for the simulation
  of cardiac electromechanics in the human left ventricle, Mathematics in
  Engineering 1 (2018) 1--37.

\bibitem{pfaller2019importance}
M.~Pfaller, J.~H{\"o}rmann, M.~Weigl, A.~Nagler, R.~Chabiniok, C.~Bertoglio,
  W.~Wall, The importance of the pericardium for cardiac biomechanics: from
  physiology to computational modeling, Biomechanics and Modeling in
  Mechanobiology 18 (2019) 503--529.

\bibitem{strocchi2020simulating}
M.~Strocchi, M.~Gsell, C.~Augustin, O.~Razeghi, C.~Roney, A.~Prassl,
  E.~Vigmond, J.~Behar, J.~Gould, C.~Rinaldi, et~al., Simulating ventricular
  systolic motion in a four-chamber heart model with spatially varying {R}obin
  boundary conditions to model the effect of the pericardium, Journal of
  {B}iomechanics 101 (2020) 109645.

\bibitem{quarteroni2019cardiovascular}
A.~Quarteroni, L.~Ded{\`e}, A.~Manzoni, C.~Vergara, Mathematical modelling of
  the human cardiovascular system: data, numerical approximation, clinical
  applications, Cambridge Monographs on Applied and Computational Mathematics,
  Cambridge University Press, 2019.

\bibitem{quarteroni2009numerical}
A.~Quarteroni, Numerical models for differential problems, Vol.~2, Springer,
  2009.

\bibitem{Africa2019}
P.~Africa, Scalable adaptive simulation of organic thin-film transistors, Ph.D.
  thesis, Politecnico di Milano (2019).

\bibitem{BursteddeWilcoxGhattas11}
B.~Carsten, C.~Lucas, G.~Omar, {\texttt{p4est}}: Scalable algorithms for
  parallel adaptive mesh refinement on forests of octrees, SIAM Journal on
  Scientific Computing 33 (2011) 1103--1133.

\bibitem{quarteroni2010numerical}
A.~Quarteroni, R.~Sacco, F.~Saleri, Numerical mathematics, Vol.~37, Springer
  Science \& Business Media, 2010.

\bibitem{niederer2011simulating}
S.~Niederer, L.~Mitchell, N.~Smith, G.~Plank, Simulating human cardiac
  electrophysiology on clinical time-scales, Frontiers in Physiology 2 (2011)
  14.

\bibitem{benzi2005numerical}
M.~Benzi, G.~Golub, J.~Liesen, et~al., Numerical solution of saddle point
  problems, Acta Numerica 14 (2005) 1--137.

\bibitem{regazzoni2021emulator}
F.~Regazzoni, A.~Quarteroni, Accelerating the convergence to a limit cycle in
  3d cardiac electromechanical simulations through a data-driven 0d emulator,
  Computers in Biology and Medicine (2021) 104641.

\bibitem{zygote2014}
Z.~M.~G. Inc., {Zygote solid 3d heart generation II developement report.},
  Technical Report (2014).

\bibitem{antiga2008image}
L.~Antiga, M.~Piccinelli, L.~Botti, B.~Ene-Iordache, A.~Remuzzi, D.~Steinman,
  An image-based modeling framework for patient-specific computational
  hemodynamics, Medical \& Biological Engineering \& Computing 46~(11) (2008)
  1097--1112.

\bibitem{fedele2021polygonal}
M.~Fedele, A.~Quarteroni, Polygonal surface processing and mesh generation
  tools for the numerical simulation of the cardiac function, International
  Journal for Numerical Methods in Biomedical Engineering 37~(4) (2021) e3435.

\bibitem{durrer1970total}
D.~Durrer, R.~Van~Dam, G.~Freud, M.~Janse, F.~Meijler, R.~Arzbaecher, Total
  excitation of the isolated human heart, Circulation 41~(6) (1970) 899--912.

\bibitem{dealII91}
D.~Arndt, W.~Bangerth, T.~Clevenger, D.~Davydov, M.~Fehling, D.~Garcia-Sanchez,
  G.~Harper, T.~Heister, L.~Heltai, M.~Kronbichler, R.~Kynch, M.~Maier, J.-P.
  Pelteret, B.~Turcksin, D.~Wells, {The \texttt{deal.II} Library, Version 9.1},
  Journal of Numerical Mathematics (2019).

\bibitem{lombaert2012human}
H.~Lombaert, J.~Peyrat, P.~Croisille, S.~Rapacchi, L.~Fanton, F.~Cheriet,
  P.~Clarysse, I.~Magnin, H.~Delingette, N.~Ayache, Human atlas of the cardiac
  fiber architecture: study on a healthy population, IEEE Transactions on
  Medical Imaging 31~(7) (2012) 1436--1447.

\bibitem{anderson2009three}
R.~Anderson, M.~Smerup, D.~Sanchez-Quintana, M.~Loukas, P.~Lunkenheimer, The
  three-dimensional arrangement of the myocytes in the ventricular walls,
  Clinical Anatomy: The Official Journal of the American Association of
  Clinical Anatomists and the British Association of Clinical Anatomists 22~(1)
  (2009) 64--76.

\end{thebibliography}
